\documentclass[final,3p,sort&compress]{elsarticle}

\usepackage{amssymb,amsmath,amsthm,amsbsy}
\usepackage{comment}
\usepackage{physics}
\usepackage{graphicx}
\usepackage{epstopdf}
\usepackage{epsfig}
\usepackage{subfigure}
\usepackage{multicol}
\usepackage{multirow}
\usepackage{setspace}
\usepackage{color}
\usepackage{tcolorbox}
\usepackage{url}
\usepackage[mathlines, running]{lineno}

\usepackage{algorithm}
\usepackage[noend]{algpseudocode}
\algnewcommand{\LineComment}[1]{\vspace{.0625in}\State \textcolor{gray}{\texttt{\# #1}}}
\algrenewcomment[1]{\hfill\textcolor{gray}{\texttt{\# #1}}}

\newcommand{\trp}{^{\mathsf{T}}}
\numberwithin{equation}{section}

\journal{Computer Methods in Applied Mechanics and Engineering}

\begin{document}

\begin{frontmatter}

\title{Bayesian operator inference for data-driven reduced-order modeling}

\author[add1]{Mengwu Guo\corref{cor1}}
\ead{m.guo@utwente.nl}
\cortext[cor1]{Corresponding author.}

\author[add2]{Shane A. McQuarrie}
\ead{shanemcq@utexas.edu}

\author[add2]{Karen E. Willcox}
\ead{kwillcox@oden.utexas.edu}

\address[add1]{Department of Applied Mathematics, University of Twente}
\address[add2]{Oden Institute for Computational Engineering and Sciences, The University of Texas at Austin}

\begin{abstract}
This work proposes a Bayesian inference method for the reduced-order modeling of time-dependent systems. Informed by the structure of the governing equations, the task of learning a reduced-order model from data is posed as a Bayesian inverse problem with Gaussian prior and likelihood. The resulting posterior distribution characterizes the operators defining the reduced-order model, hence the predictions subsequently issued by the reduced-order model are endowed with uncertainty. The statistical moments of these predictions are estimated via a Monte Carlo sampling of the posterior distribution. Since the reduced models are fast to solve, this sampling is computationally efficient. Furthermore, the proposed Bayesian framework provides a statistical interpretation of the regularization term that is present in the deterministic operator inference problem, and the empirical Bayes approach of maximum marginal likelihood suggests a selection algorithm for the regularization hyperparameters. The proposed method is demonstrated on two examples: the compressible Euler equations with noise-corrupted observations, and a single-injector combustion process.
\end{abstract}

\end{frontmatter}

\section{Introduction}

\noindent
The objective of this work is to enable uncertainty quantification for data-driven reduced-order modeling. A reduced-order model (ROM) is a low-dimensional surrogate for a high-dimensional physical system, such as the many natural and engineering systems that are governed by partial differential equations. Ideally, the computational cost of a ROM should be reduced compared to that of a high-fidelity model of the original system, yet without compromising the solution accuracy. Model reduction strategies typically construct a ROM by identifying a predominant, low-dimensional subspace that well approximates the solution manifold, then representing the state dynamics in the subspace \cite{BGW2015pmorSurvery,quarteroni2014reduced}. Classical reduced-order models are constructed by explicitly projecting high-fidelity operators onto the reduced subspace. Such models are called \emph{intrusive} because the high-fidelity operators must be extracted from a numerical solution code. In contrast, \emph{non-intrusive} data-driven ROMs learn low-dimensional representations of systems primarily from solution data, and do not require intrusive access to a high-fidelity simulator. Non-intrusive ROMs are advantageous for applications with well-established legacy codes that can be readily executed but not easily modified. On the other hand, the extensive use of data introduces uncertainties beyond those originating from the modeling decisions of low-dimensional approximation and makes it difficult to obtain classical \emph{a posteriori} error estimates as in \cite{hesthaven2016certified}. In particular, measurement noise in the data and misspecification of the model structure may impact the accuracy and robustness of non-intrusive ROMs. It is therefore crucial that validation efforts for data-driven ROMs quantify the errors and uncertainties associated with non-intrusive learning.

There are several data-driven methods for constructing ROMs, including dynamic mode decomposition \cite{schmid2010dynamic}, operator inference to learn a ROM with polynomial structure inspired by known governing equations \cite{GW2021learning,PW2016operatorInference}, sparse regression to identify reduced-order latent dynamics \cite{champion2019sindy,schaeffer2017learning}, manifold learning using deep auto-encoders \cite{Lee2020autoencoder}, data-driven approximation of time-integration schemes \cite{zhuang2021model}, and Gaussian process modeling for parametric reduced basis methods \cite{guo2019data}. This paper focuses on operator inference (OpInf), a non-intrusive, projection-based model reduction technique introduced in \cite{PW2016operatorInference}. OpInf learns a ROM by inferring reduced-order operators through a deterministic linear regression. The regression requires solution data of the states and their corresponding time derivatives; if needed, the time derivatives are typically estimated from the states with finite differences. Importantly, numerical regularization is often required in the inference problem to address ill-conditioning, over-fitting, and/or model inadequacy \cite{MHW2021regOpInfCombustion,SKHW2020romCombustion}, but the regularization selection remains a subject of further investigation. Errors introduced by time derivative approximation, noise in the solution data, and/or numerical regularization have not been directly addressed in the literature, with the exception of \cite{uy2021activeOpInf}, which considers an active sampling strategy for overcoming noise pollution in the training data. Even without noise in the data, establishing error bounds or uncertainty guarantees in OpInf is challenging because the learning is driven by data, although \cite{uy2021probabilistic} establishes probabilistic error estimates for linear OpInf ROMs in the case where the reduced-order operators can be ensured to be the same as those determined by explicit projection as in \cite{Peherstorfer2020reprojection}.

The main contribution of this work is to establish a probabilistic formulation of OpInf, thereby facilitating uncertainty quantification for non-intrusive ROMs of time-dependent systems. Our strategy accounts for approximation errors and data noise, endows the resulting ROM with uncertainty bounds, and provides statistical insight into the numerical regularization. Specifically, the task of inferring reduced-order operators is posed as a Bayesian inversion problem \cite{box2011bayesian}, which incorporates prior knowledge and observational information to update the probabilistic distribution of unknown system parameters. Bayesian approaches are especially useful in data-driven modeling and analysis \cite{chkrebtii2016bayesian} and, in particular, have recently shown promise in estimating the operators in differential equations, see, e.g., \cite{HBDK2021uqsindy,pan2015sparse,zhang2018robust}. In our case, multivariate normal distributions are adopted to define the prior and likelihood when Bayes' rule is applied, yielding a Bayesian linear regression to estimate the reduced-order operators. Different from the deterministic OpInf workflow that gives point estimates only, the proposed Bayesian version derives posterior distributions of the learned reduced-order operators, whose uncertainties are subsequently propagated to solutions through time integration. In other words, endowing the reduced-order operators with probabilistic descriptions results in ROM predictions with quantifiable uncertainties. Additionally, our probabilistic treatment of the regression leads to a selection algorithm for the regularization hyperparameters.

We present this work in honor of the lifetime academic achievements of Professor J. Tinsley Oden, especially his seminal work in highlighting the importance of \emph{verification \& validation} \cite{oden2011fembook,oden2004verification,oden2002estimation} and \emph{uncertainty quantification} \cite{oden2017predictive,oden2005uncertainty} to achieve predictive models in computational engineering and sciences \cite{oden2018predictivemodelling}. He laid the mathematical foundation of predictive computational science \cite{oden2017predictive} using Bayesian methods and has utilized the framework in a substantial range of scientific, medical and engineering applications \cite{oden2015bayesatomistic,oden2013bayestumor,oden2020bayescovid,oden2015bayesdamage}. Professor Oden's pioneering, groundbreaking research on predictive computational science sets the stage for this work and inspires us to enable the uncertainty quantification and model validation of data-driven ROMs through a Bayesian approach.

The paper is structured as follows. Section~\ref{sec:opinf-intro} introduces OpInf in its established deterministic framework. Section~\ref{sec:opinf-bayes} endows OpInf with a Bayesian perspective by establishing a posterior distribution for the operators of a reduced-order system, and Section~\ref{sec:regularization} discusses the implications for the numerical regularization that is needed in the ROM learning formulation. Numerical results for two examples are presented in Section~\ref{sec:results}, and Section~\ref{sec:conclusion} makes concluding remarks.

\section{Deterministic operator inference}
\label{sec:opinf-intro}

We consider problems governed by partial differential equations, written after spatial discretization in the form
\begin{align}
    \label{eq:fullorder}
    \frac{\textup{d}}{\textup{d}t}{\vb{q}}(t)
    = {\vb{F}}({\vb{q}}(t),{\vb{u}}(t))\,,
    \quad
    \vb{q}(t_0) = \vb{q}_0\,,
    \quad
    t\in [t_0, t_f]\,,
\end{align}
in which $\vb{q}(t)\in \mathbb{R}^n$ is the state vector collecting values of $n$ degrees of freedom at time $t$,
$\vb{u}(t)\in \mathbb{R}^m$ collects $m$ input values at time $t$ (e.g., encoding boundary conditions or time-dependent source terms),
$\vb{F}:\mathbb{R}^{n}\times\mathbb{R}^{m}\to\mathbb{R}^{n}$ maps state/input pairs to the state evolution,
$t_0$ and $t_f$ are the initial and final time, respectively, and
$\vb{q}_0$ is the prescribed initial condition.
A projection-based ROM approximates the full-order system~\eqref{eq:fullorder} in a subspace with reduced dimensionality $r\ll n$, spanned by a collection of $r$ basis vectors $\vb{V}=[\vb{v}_1~\vb{v}_2~\cdots~\vb{v}_r]\in\mathbb{R}^{n\times r}$.
One widely used approach to construct the basis $\vb{V}$ is the proper orthogonal decomposition (POD) \cite{Berkooz1993, GPT1999vortexsheddingPOD, sirovich1987turbulence}.
For a collection of known state snapshots $\vb{q}_{j} := \vb{q}(t_j)$ and inputs $\vb{u}_{j} := \vb{u}(t_{j})$, $j=0,1,\cdots,k-1$, at $k$ time instances $t_0< t_1 <\cdots<t_{k-1}$, organized as $\vb{Q} =  [~\vb{q}_{0}~\vb{q}_{1}~\cdots~\vb{q}_{k-1}~]\in\mathbb{R}^{n\times k}$ and $\vb{U} =  [~\vb{u}_{0}~\vb{u}_{1}~\cdots~\vb{u}_{k-1}~]\in\mathbb{R}^{m\times k}$,
POD takes advantage of the singular value decomposition $\vb{Q}=\vb*{\Phi}\vb*{\Xi}\vb*{\Psi}\trp$ by using the first $r$ columns of $\vb*{\Phi}$ as the $r$-dimensional reduced basis, i.e., $\vb{V}=\vb*{\Phi}_{:,1:r}$.
Projection-based methods, such as those based on POD, have the key advantage of producing ROMs that inherit structure from the full-order system~\eqref{eq:fullorder} \cite{GPT1999vortexsheddingPOD}.
This observation motivates us to consider, for a wide class of problems, ROMs with the following polynomial structure:
\begin{equation}\label{eq:reducedorder}
\frac{\textup{d}}{\textup{d}t} \hat{\vb{q}}(t)= \hat{\vb{A}}\hat{\vb{q}}(t) + \hat{\vb{H}}[\hat{\vb{q}}(t)\otimes \hat{\vb{q}}(t) ]+\hat{\vb{B}}{\vb{u}}(t) + \hat{\vb{c}}\,,\quad \hat{\vb{q}}(t_0) = \vb{V}\trp\vb{q}_0\,, \quad t\in [t_0, t_f]\,,
\end{equation}
where
$\hat{\vb{q}}(t) \in \mathbb{R}^{r}$ is the reduced state,
$\otimes$ denotes the Kronecker product \cite{vanLoan2000kronecker}, and
$\hat{\vb{A}}\in \mathbb{R}^{r\times r}$,
$\hat{\vb{H}}\in \mathbb{R}^{r\times r^2}$,
$\hat{\vb{B}}\in\mathbb{R}^{r\times m}$, and
$\hat{\vb{c}}\in\mathbb{R}^{r}$ are reduced-order operators that together comprise a low-dimensional representation of ${\vb{F}}$.
The full- and reduced-order states are related by a linear approximation $\mathbf{q}(t) \approx \mathbf{V}\hat{\mathbf{q}}(t)$.

The data-driven \emph{operator inference} (OpInf) method \cite{PW2016operatorInference} employs the projection-based setting of \eqref{eq:reducedorder}, but utilizes a regression problem to find the reduced-order operators that best match the snapshot data in a minimum-residual sense.
For non-polynomial problems, training data may be pre-processed with variable transformations to induce an exact or approximate polynomial structure as needed \cite{QKPW2020liftAndLearn}.
Such a non-intrusive approach avoids entailing intrusive queries and access to source code. The regression problem to infer the reduced-order operators is written as
\begin{equation}\label{eq:opinf}
    \min_{\hat{\vb{A}},\hat{\vb{H}},\hat{\vb{B}},\hat{\vb{c}}}
    \left\{\sum_{j=0}^{k-1}\left\|
        \hat{\vb{A}}\hat{\vb{q}}_{j}
        + \hat{\vb{H}}[\hat{\vb{q}}_{j} \otimes \hat{\vb{q}}_{j}]
        + \hat{\vb{B}}{\vb{u}}_{j}
        + \hat{\vb{c}}
        - \dot{\hat{\vb{q}}}_{j}
    \right\|_2^2
    + \mathcal{P}\left([~\hat{\vb{A}}~~\hat{\vb{H}}~~\hat{\vb{B}}~~\hat{\vb{c}}~]\right)\right\}
    \,,
\end{equation}
in which the first term is a loss function defined by the residual of~\eqref{eq:reducedorder} at the $k$ time instances $t_0< t_1 <\cdots<t_{k-1}$, and the second term $\mathcal{P}$ represents a regularization. Here the reduced-state dataset $\{\hat{\vb{q}}_{j}\}_{j=0}^{k-1}$ is constructed by evaluating the projection of known snapshots onto the reduced basis $\vb{V}$, i.e., $\hat{\vb{q}}_{j} = \vb{V}\trp\vb{q}_{j}$, and the time-derivative data $\{\dot{\hat{\vb{q}}}_{j}\}_{j=0}^{k-1}$ are then typically approximated using finite differences.
The regularization term $\mathcal{P}$ penalizes the entries of the learned quantities, which is often necessary to induce stability in the resulting ROM \cite{MHW2021regOpInfCombustion}.
The training time instances are contained in $[t_0, t_{k-1}]$, a subset of the full time interval of interest $[t_0, t_f]$, and the reduced-order solution over $[t_{k},t_f]$ will be predictive through the inferred ROM.

\vspace{3mm}\noindent\textbf{Remark 1:}
The quadratic interaction $\hat{\mathbf{q}}\otimes \hat{\mathbf{q}} = \text{vec}(\hat{\mathbf{q}}\hat{\mathbf{q}}\trp)\in \mathbb{R}^{r^2}$
has only $r(r+1)/2 < r^2$ degrees of freedom due to symmetry, i.e., the upper-triangular elements of $\hat{\mathbf{q}}\hat{\mathbf{q}}\trp$ completely prescribe $\hat{\mathbf{q}}\otimes\hat{\mathbf{q}}$.
Hence, in practice we may learn an $r \times r(r+1)/2$ quadratic operator $\hat{\mathbf{H}}$ as opposed to an $r \times r^2$ quadratic operator.
\vspace{3mm}

The least-squares problem~\eqref{eq:opinf} can be rewritten in a matrix form and decoupled into $r$ least-squares problems, one for each ordinary differential equation in the reduced system~\eqref{eq:reducedorder}, i.e.,
\begin{equation}
    \label{eq:leastsqr}
    \min_{\hat{\vb{O}}} \left\{\left\|
        \vb{D}\hat{\vb{O}}\trp - \vb{R}\trp
    \right\|_F^2
    + \mathcal{P}(\hat{\mathbf{O}}) \right\}
    =  \sum_{i=1}^{r} \min_{\hat{\vb{o}}_i}\left\{\left\|
        \vb{D}\hat{\vb{o}}_i -\vb{r}_i
    \right\|_2^2
    + \mathcal{P}_{i}(\hat{\mathbf{o}}_{i})\right\}\,,
\end{equation}
where
\begin{align}
\label{eq:matrix-definitions}
\begin{aligned}
    \hat{\vb{O}} & := [~\hat{\vb{A}}~~\hat{\vb{H}}~~\hat{\vb{B}}~~\hat{\vb{c}}~] := [~\hat{\vb{o}}_1~~\cdots~~\hat{\vb{o}}_r~]\trp\in \mathbb{R}^{r\times d(r,m)}\,,
    \\
    \vb{D} & := [~\hat{\vb{Q}}\trp~~ (\hat{\vb{Q}}\odot \hat{\vb{Q}})\trp~~\vb{U}\trp~~ \vb{1}_k~]
    \in \mathbb{R}^{k\times d(r,m)}\,,
    \\
    \hat{\vb{Q}} & := [~\hat{\vb{q}}_{0}~~\hat{\vb{q}}_{1}~~\cdots~~\hat{\vb{q}}_{k-1}~] \in \mathbb{R}^{r\times k}\,,
    \\
    \vb{R} & := [~\dot{\hat{\vb{q}}}_{0}~~\dot{\hat{\vb{q}}}_{1}~~\cdots~~\dot{\hat{\vb{q}}}_{k-1}~] := [~\vb{r}_1~~\cdots~~\vb{r}_r~]\trp\in \mathbb{R}^{r\times k}\,,
    \quad \text{and}
    \\
    \vb{U} & := [~\vb{u}_{0}~~\vb{u}_{1}~~\cdots~~\vb{u}_{k-1}~] \in \mathbb{R}^{m\times k}\,,
\end{aligned}
\end{align}
in which $d(r,m) = r + r(r+1)/2 + m + 1$,
$\vb{1}_k\in \mathbb{R}^k$ is a length-$k$ vector with all entries set to 1,
$\odot$ is the Khatri-Rao product (the column-wise Kronecker product \cite{KhatriRoa1968product,kolda2009tensor}), and
$\hat{\vb{o}}_i\in \mathbb{R}^{d(r,m)}$ and $\vb{r}_i\in\mathbb{R}^k$ are the $i$-th rows of $\hat{\vb{O}}$ and $\vb{R}$, respectively.
Here an individual regularization $\mathcal{P}_i$ is defined for each $\hat{\vb{o}}_i$, i.e., $\mathcal{P}(\hat{\vb{O}})=\sum_{i=1}^{r}\mathcal{P}_i(\hat{\vb{o}}_i)$. Thus the OpInf comprises of $r$ independent regressions with respect to the unknown reduced-order operators $\hat{\vb{o}}_i$. Note that while \eqref{eq:leastsqr} is a linear regression, the corresponding ROM \eqref{eq:reducedorder} defined by $\hat{\vb{o}}_{i}$ is a nonlinear dynamical system.

We usually have $k>d(r,m)$; in that case, and with $\mathcal{P}_{i}(\hat{\vb{o}}_{i}) = 0$, the point estimate $\tilde{\vb{o}}_{i}$ of $\hat{\vb{o}}_i$ through the least squares~\eqref{eq:leastsqr} can be yielded in an explicit form as
\begin{align}\label{eq:pointest}
    \tilde{\vb{o}}_i :=  \vb{D}^+ \vb{r}_i = (\vb{D}\trp\vb{D})^{-1}\vb{D}\trp\vb{r}_i\,,
\end{align}
$\vb{D}^+$ denoting the Moore–Penrose inverse of the matrix $\vb{D}$.
Here we assume that $\mathbf{D}$ has full column rank; the cases with ill-conditioned $\mathbf{D}\trp\mathbf{D}$ will be discussed in Section \ref{sec:regularization}.

\section{Bayesian operator inference}
\label{sec:opinf-bayes}

In this section, the regression for learning the reduced-order operators $\hat{\mathbf{O}} = [~\hat{\vb{A}}~~\hat{\vb{H}}~~\hat{\vb{B}}~~\hat{\vb{c}}~]$ is posed as a Bayesian inference problem \cite{box2011bayesian,chkrebtii2016bayesian}. We define a probabilistic ROM through the resulting posterior distribution of $\hat{\mathbf{O}}$ substituted into the form~\eqref{eq:reducedorder}, which enables uncertainty estimates through Monte-Carlo sampling. We abbreviate the proposed \emph{Bayesian operator inference} method as BayesOpInf.

\subsection{Bayesian inference of reduced-order operators}

The quadratic structure of the ROM~\eqref{eq:reducedorder} introduces a model misspecification error whenever the full-order system~\eqref{eq:fullorder} includes non-quadratic dynamics. Additionally, the quality of the least-squares estimation~\eqref{eq:leastsqr} may suffer from approximation error in the estimated time derivatives of the training snapshots and/or noise in the snapshots themselves. To account for these uncertainties, we use Bayesian inference to formulate the estimators of the reduced-order operators $\hat{\vb{o}}_i$ as probabilistic distributions,  $i = 1, \ldots, r$.

For the definition of a likelihood distribution, we model the residual error in the data-driven scheme~\eqref{eq:leastsqr} as independent Gaussian noise, i.e., we consider
\begin{align}
    \label{eq:leastsqr-bayes}
    \vb{r}_{i} = \vb{D}\hat{\vb{o}}_{i} + \boldsymbol{\epsilon}_{i},
    \qquad
    p(\boldsymbol{\epsilon}_i)
    = \mathcal{N}(\boldsymbol{\epsilon}_i|\vb{0}_{k},\sigma_i^2\vb{I}_{k})\,,
\end{align}
in which $\vb{0}_{k}$ is the $k$-dimensional zero vector and $\vb{I}_{k}$ denotes the $k \times k$ identity matrix. Here the noise $\vb*{\epsilon}_i\in\mathbb{R}^{k}$ is a random vector consisting of the error terms in the observations at $k$ time instances. All these noise measurements are assumed to follow independent normal distributions with zero mean and variance $\sigma_i^2$. Hence the likelihood distribution for the training data pair $(\vb{D},\vb{r}_{i})$ is written as
\begin{align}\label{eq:likelihood}
    p(\vb{r}_i|~\vb{D},\hat{\vb{o}}_i,\sigma_i^2)
    = \mathcal{N}(\vb{r}_i |~\vb{D}\hat{\vb{o}}_i,\sigma_i^2\vb{I}_{k})\,.
\end{align}
The choice to model the entries of $\boldsymbol{\epsilon}_{i}$ as independent random variables (i.e., assuming the noise is not correlated across time) is convenient but not a requirement. An alternative, kernel-based model of the residual error is provided in Appendix~1, which accounts for correlations over time and constructs a surrogate model of the ROM closure.

We also define the prior distribution of $\hat{\vb{o}}_i$ as a multivariate normal distribution,
\begin{align}\label{eq:prior}
    p(\hat{\vb{o}}_i|~\sigma_i^2, \vb*{\lambda}_i)
    = \mathcal{N}(\hat{\vb{o}}_i|~\vb*{\beta}_i,\sigma_i^2\mathrm{diag}(\vb*{\lambda}_i)^{-1})\,,
\end{align}
where $\vb*{\beta}_i\in \mathbb{R}^{d(r,m)}$ is a pre-defined mean vector and $\vb*{\lambda}_i \in \mathbb{R}^{d(r,m)}$ is a vector of positive numbers hyperparametrizing the prior variances of $\hat{\vb{o}}_i$. Note that the prior covariance is a diagonal matrix, implying that the prior does not impose correlations between the components of $\hat{\vb{o}}_i$.

Conditioning on the training data, Bayes' rule gives the posterior distribution $p(\hat{\vb{o}}_i | ~\vb{D},\vb{r}_i,\sigma_i^2, \vb*{\lambda}_i)$ as follows:
\begin{equation}
\begin{split}\label{eq:post}
p(\hat{\vb{o}}_i | ~\vb{D},\vb{r}_i,\sigma_i^2, \vb*{\lambda}_i) & \propto p(\vb{r}_i|~\vb{D},\hat{\vb{o}}_i,\sigma_i^2)~ p(\hat{\vb{o}}_i|~\sigma_i^2,\vb*{\lambda}_i)
\\
& = \mathcal{N}(\vb{r}_i |~\vb{D}\hat{\vb{o}}_i,\sigma_i^2\vb{I}_k) ~\mathcal{N}(\hat{\vb{o}}_i|~\vb*{\beta}_i,\sigma_i^2\mathrm{diag}(\vb*{\lambda}_i)^{-1})
\\
& = \mathcal{N}(\vb{D}\hat{\vb{o}}_i |~\vb{r}_i,\sigma_i^2\vb{I}_k) ~\mathcal{N}(\hat{\vb{o}}_i|~\vb*{\beta}_i,\sigma_i^2\mathrm{diag}(\vb*{\lambda}_i)^{-1})
\\
& = \mathcal{N}\left(\vb{r}_i |~ \vb{D}\vb*{\beta}_i, \sigma_i^2\vb{D}\mathrm{diag}(\vb*{\lambda}_i)^{-1}\vb{D}\trp+\sigma_i^2\vb{I}_k\right)~\mathcal{N}(\hat{\vb{o}}_i |~ \vb*{\mu}_i,\vb*{\Sigma}_i)\,,\quad \text{i.e.,}
\\
p(\hat{\vb{o}}_i | ~\vb{D},\vb{r}_i,\sigma_i^2, \vb*{\lambda}_i) & = \mathcal{N}(\hat{\vb{o}}_i | ~\vb*{\mu}_i,\vb*{\Sigma}_i)\,,
\end{split}
\end{equation}
in which
\begin{equation}
    \begin{split}\label{eq:musigma}
    \vb*{\Sigma}_i
    & = \vb*{\Sigma}_i(\vb{D},\sigma_i^2, \vb*{\lambda}_i)
    = \sigma_i^2\left[ \mathrm{diag}(\vb*{\lambda}_i) +\vb{D}\trp\vb{D} \right]^{-1},
    \\
    \vb*{\mu}_i
    & = \vb*{\mu}_i(\vb{r}_i,\vb{D},\vb*{\lambda}_i)
    = \left[ \mathrm{diag}(\vb*{\lambda}_i) +\vb{D}\trp\vb{D} \right]^{-1} \left( \mathrm{diag}(\vb*{\lambda}_i)\vb*{\beta}_i+\vb{D}\trp\vb{r}_i  \right)
    \\
    & = \vb*{\beta}_i + \underbrace{\left[ \mathrm{diag}(\vb*{\lambda}_i) +\vb{D}\trp\vb{D} \right]^{-1}\vb{D}\trp(\vb{r}_i -\vb{D}\vb*{\beta}_i)}_{\delta\vb*{\mu}_i(\vb{r}_i,\vb{D},\vb*{\lambda}_i)}\,,
    \end{split}
\end{equation}
all consistent with the results of Bayesian linear regression \cite{box2011bayesian}.

To marginalize $(\sigma_i^2,\vb*{\lambda}_i)$ in the distribution $p(\hat{\vb{o}}_i | ~\vb{D},\vb{r}_i,\sigma_i^2, \vb*{\lambda}_i)$ and obtain $p(\hat{\vb{o}}_i | ~\vb{D},\vb{r}_i)$ as the inference of $\hat{\vb{o}}_i$, we adopt an empirical Bayes method---maximum log marginal likelihood  \cite{rasmussen2006gaussian}---to estimate a set of optimal values of $(\sigma_i^2,\vb*{\lambda}_i)$. Mathematically, the maximization problem is written as
\begin{equation}\label{eq:maxlogmlhd}
\begin{split}
(\sigma^{*2}_i,\vb*{\lambda}_i^*) =  \underset{\sigma_i^2,\vb*{\lambda}_i}{\textrm{arg\,max}}&~\log~p(\vb{r}_i|~\vb{D},\sigma_i^2,\vb*{\lambda}_i)
\\
= \underset{\sigma_i^2,\vb*{\lambda}_i}{\textrm{arg\,max}}&~\log\int p(\vb{r}_i|~\vb{D},\hat{\vb{o}}_i,\sigma_i^2)~ p(\hat{\vb{o}}_i|~\sigma_i^2,\vb*{\lambda}_i)~\dd\hat{\vb{o}}_i \\
 \overset{\text{from~\eqref{eq:post}}}{=}  \underset{\sigma_i^2,\vb*{\lambda}_i}{\textrm{arg\,max}}&~\log~\mathcal{N}\left(\vb{r}_i | ~\vb{D}\vb*{\beta}_i, \sigma_i^2\left(\vb{D}\mathrm{diag}(\vb*{\lambda}_i)^{-1}\vb{D}\trp+\vb{I}_k\right)\right)
\\
= \underset{\sigma_i^2,\vb*{\lambda}_i}{\textrm{arg\,max}}&~ -\frac{1}{2\sigma_i^2}\|\vb{r}_i-\vb{D}\vb*{\mu}_i(\vb{r}_i,\vb{D},\vb*{\lambda}_i) \|_2^2-\frac{1}{2\sigma_i^2} \left\|\mathrm{diag}\left(\vb*{\lambda}_i\right)^{1/2}\delta\vb*{\mu}_i(\vb{r}_i,\vb{D},\vb*{\lambda}_i)\right\|_2^2\\
& ~ -\frac{1}{2}\log~\left|\mathrm{diag}(\vb*{\lambda}_i)+\vb{D}\trp\vb{D}\right| +\frac{1}{2}\log(\vb*{\lambda}_i)\trp\vb{1}_{d(r,m)} -\frac{k}{2}\log(\sigma_i^2)-\frac{k}{2}\log(2\pi)\,,
\end{split}
\end{equation}
where $|\cdot|$ indicates the matrix determinant. Thus the final posterior distribution $p(\hat{\vb{o}}_i | ~\vb{D},\vb{r}_i)$  is approximated as
\begin{equation}
\label{eq:final-operator-posterior}
p(\hat{\vb{o}}_i | ~\vb{D},\vb{r}_i) = p(\hat{\vb{o}}_i | ~\vb{D},\vb{r}_i,\sigma_i^{*2}, \vb*{\lambda}^*_i) =\mathcal{N}\left( \hat{\vb{o}}_i|~  \vb*{\mu}_i(\vb{r}_i,\vb{D},\vb*{\lambda}^*_i), \vb*{\Sigma}_i(\vb{D},\sigma_i^{*2}, \vb*{\lambda}^*_i) \right)\,.
\end{equation}
Furthermore, by taking partial derivative of the objective function in \eqref{eq:maxlogmlhd} with respect to $\sigma_i^2$ and setting it to zero, we see that the maximum is achieved at
\begin{equation}\label{eq:noise}
    \sigma_i^{*2} = \frac{\|\vb{r}_i-\vb{D}\vb*{\mu}_i\|_2^2+\left\|\mathrm{diag}\left(\vb*{\lambda}_i^*\right)^{1/2}\delta\vb*{\mu}_i\right\|_2^2}{k}\,.
\end{equation}

\vspace{3mm}
\noindent\textbf{Remark 2:} A pragmatic choice is to set $\vb*{\beta}_i = \vb{0}$ and thus $\delta\vb*{\mu}_i = \vb*{\mu}_i $. When we assume that the prior variance values are very large ($\vb*{\lambda}_i\to \vb{0}$), i.e., the prior of $\hat{\vb{o}}_i$ is uninformative, we have $\vb*{\mu}_i = \tilde{\vb{o}}_i = \vb{D}^+ \vb{r}_i $, $\vb*{\Sigma}_i = \sigma_i^{*2} (\vb{D}\trp\vb{D})^{-1}$, and $\sigma_i^{*2} = \|\vb{r}_i-\vb{D}\tilde{\vb{o}}_i \|_2^2/k$. In this case, the posterior mean vector of $\hat{\vb{o}}_i$ coincides with the original least squares estimate $\tilde{\vb{o}}_i $ as in deterministic OpInf, and $\sigma_i^{*2}$ is the mean squared error of the linear regression.

\vspace{3mm}
\noindent\textbf{Remark 3:} If the pre-defined prior mean vector of $\hat{\vb{o}}_i$ is taken as its point estimate $\tilde{\vb{o}}_i$, i.e., $\vb*{\beta}_i = \tilde{\vb{o}}_i = \vb{D}^+ \vb{r}_i $, then we have $\delta \vb*{\mu}_i=\vb{0}$ and $\vb*{\mu}_i=\tilde{\vb{o}}_i $, as well as $\sigma_i^{*2} = \|\vb{r}_i-\vb{D}\tilde{\vb{o}}_i \|_2^2/k$. Moreover, by simplifying \eqref{eq:maxlogmlhd} the estimate of $\vb*{\lambda}_i$ is given by
\begin{equation}
\begin{split}
\vb*{\lambda}_i^*  & =  \underset{\vb*{\lambda}_i}{\textrm{arg\,max}}~
- \frac{1}{2} \log~ |\mathrm{diag}(\vb*{\lambda}_i)+\vb{D}\trp\vb{D}|+\frac{1}{2}\log ~ |\mathrm{diag}(\vb*{\lambda}_i)| \\
& = \underset{\vb*{\lambda}_i}{\textrm{arg\,max}}~
- \frac{1}{2} \log~ \left|\vb{I}+\mathrm{diag}(\vb*{\lambda}_i)^{-\frac{1}{2}}\vb{D}\trp\vb{D}\mathrm{diag}(\vb*{\lambda}_i)^{-\frac{1}{2}} \right| \,.
\end{split}
\end{equation}
Since $\mathrm{diag}(\vb*{\lambda}_i)^{-\frac{1}{2}}\vb{D}\trp\vb{D}\mathrm{diag}(\vb*{\lambda}_i)^{-\frac{1}{2}}$ is semi-positive definite, $\inf_{\vb*{\lambda}_{i}} \left|\vb{I}+\mathrm{diag}(\vb*{\lambda}_i)^{-\frac{1}{2}}\vb{D}\trp\vb{D}\mathrm{diag}(\vb*{\lambda}_i)^{-\frac{1}{2}} \right|= 1$, which is achieved as $\mathrm{diag}(\vb*{\lambda}_i)^{-1} \to \vb{0}$ ($\vb*{\lambda}_i \to \infty^{d(r,m)}$), i.e., when the prior of $\hat{\vb{o}}_i$ is deterministic. In this case, the Bayesian formulation reduces to deterministic OpInf, indicating that the deterministic version essentially forms an optimal solution, provided that the Gram matrix $\vb{D}\trp\vb{D}$ is well-conditioned.

\subsection{A probabilistic reduced-order model}

\begin{figure}
    \centering
    \includegraphics[width=.9\textwidth]{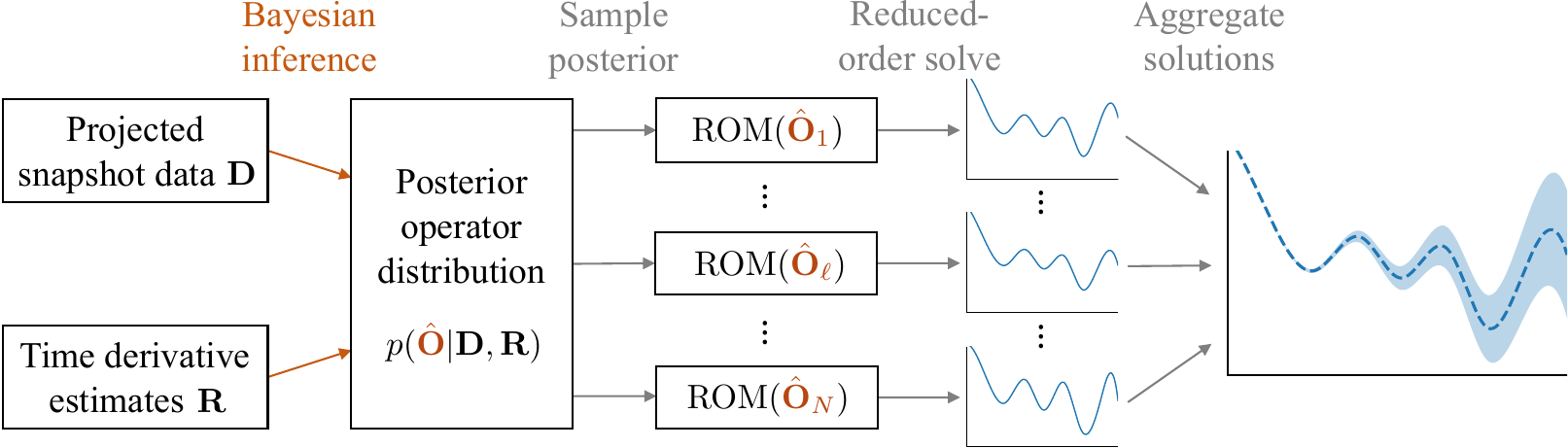}
    \caption{The BayesOpInf workflow. To derive the probabilistic ROM, infer the multivariate normal posterior distribution for the operator matrix $\hat{\mathbf{O}}$ from snapshot and time derivative data. This is a one-time, offline cost. To evaluate the probabilistic ROM in the online stage, draw samples from the posterior distribution of $\hat{\mathbf{O}}$ and solve the corresponding ROM realizations, then calculate the statistics of these ROM solutions.
    }
    \label{fig:bayesopinf-workflow}
\end{figure}

Because the OpInf regression naturally decouples along the rows of the reduced operator matrix $\hat{\vb{O}} =  [~\hat{\vb{A}}~~\hat{\vb{H}}~~\hat{\vb{B}}~~\hat{\vb{c}}~]$, we likewise treat the Bayesian inference of each row $\hat{\vb{o}}_i$ independently, $1\leq i \leq r$. Hence $\hat{\vb{O}}$ is described by the posterior distribution
\begin{equation}
p(\hat{\vb{O}}|~\vb{D},\vb{R}) = \prod_{i=1}^{r} p(\hat{\vb{o}}_i|~ \vb{D},\vb{r}_i)\,.
\end{equation}
Given $\hat{\mathbf{O}}$, we solve for the reduced state $\hat{\vb{q}}(t)$ through the corresponding reduced-order system~\eqref{eq:reducedorder} over $[t_0,t_f]$. The reduced state can therefore be viewed as a stochastic process because it depends on the random variable $\hat{\vb{O}}$, with distribution
\begin{equation}
    p(\hat{\vb{q}}(t)|~\vb{D}, \vb{R})
    = \int p(\hat{\vb{q}}(t)|~\hat{\vb{O}})p(\hat{\vb{O}}|~\vb{D},\vb{R})~\dd\hat{\vb{O}}
    = \int  p(\hat{\vb{q}}(t)|~\hat{\vb{o}}_1,\ldots,\hat{\vb{o}}_r) \prod_{i=1}^{r} p(\hat{\vb{o}}_i|~ \vb{D},\vb{r}_i)~\dd\hat{\vb{o}}_i\,,
\end{equation}
in which all the $\hat{\vb{o}}_i$'s are marginalized. Since \eqref{eq:reducedorder} is computationally inexpensive to integrate in time, we use Monte Carlo sampling over the operator posterior $p(\hat{\vb{O}}|~\vb{D},\vb{R})$ to estimate the mean function and second-order moments of $\hat{\vb{q}}(t)$ (see Figure~\ref{fig:bayesopinf-workflow}). In fact, $p(\hat{\vb{q}}(t)|\,\hat{\vb{O}})$ represents a deterministic model given by the reduced operators $\hat{\vb{O}}$, hence, one can write the reduced state solution $\hat{\vb{q}}(t)$ determined by the reduced operators $\hat{\vb{O}}$ as a function of $\hat{\vb{O}}$, i.e., $\hat{\vb{q}}(t;\hat{\vb{O}})$. Note that the solution mean and the solution determined by the posterior operator mean are not necessarily equal to each other, i.e., $\mathbb{E}[\hat{\vb{q}}(t;\hat{\vb{O}}|~\vb{D},\vb{R})]\neq \hat{\vb{q}}(t;\mathbb{E}[\hat{\vb{O}}|~\vb{D},\vb{R}])$.
We will show, however, that in practice the discrepancy between these two quantities is often small, even in complex scenarios.

Algorithm~\ref{alg:BayesOpInf} details the BayesOpInf procedure with $\boldsymbol{\beta}_{i} = \mathbf{0}$, $\delta\boldsymbol{\mu}_{i} = \boldsymbol{\mu}_{i}$ (see Remark~2), and a fixed choice of the prior variance vectors $\boldsymbol{\lambda}_{i}$, $i=1,\ldots,r$. The mean vectors $\boldsymbol{\mu}_{i}$ and covariance matrices $\boldsymbol{\Sigma}_{i}$ are determined via \eqref{eq:musigma} and \eqref{eq:noise}. We discuss the selection of $\boldsymbol{\lambda}_{i}$ in Section~\ref{sec:regularization}.

\begin{algorithm}
\begin{algorithmic}[1]
\Procedure{BayesOpInf}{
    Data matrix $\mathbf{D}\in\mathbb{R}^{k\times d(r,m)}$,
    projected time derivatives $\mathbf{R} = [~\mathbf{r}_{1}~\cdots~\mathbf{r}_{r}~]\trp\in\mathbb{R}^{r\times k}$,
    regularization parameters $\boldsymbol{\lambda}_{1},\ldots,\boldsymbol{\lambda}_{r}\in\mathbb{R}^{d(r,m)}$,
    projected initial condition $\hat{\mathbf{q}}_{0}\in\mathbb{R}^{r}$,
    final time $t_{f} > t_{0}$,
    number of Monte Carlo samples $N\in\mathbb{N}$
    }

    \LineComment{Offline phase:~construct operator posterior $p(\hat{\mathbf{O}}|\:\mathbf{D},\mathbf{R})$.}
    \For{$i= 1, \ldots, r$}
        \State $\boldsymbol{\mu}_{i} \gets \arg\min_{\boldsymbol{\mu}_{i}}\left\|\mathbf{D}\boldsymbol{\mu}_{i} - \mathbf{r}_{i}\right\|_{2}^{2} + \|\text{diag}(\boldsymbol{\lambda}_{i})^{1/2}\boldsymbol{\mu}_{i}\|_{2}^{2}$
            \label{step:bayesopinf-mean}
            \Comment{Mean of operator row $\hat{\mathbf{o}}_{i}$.}
        \State $\sigma_{i}^{*2} \gets \frac{1}{k}\left(\left\|\mathbf{D}\boldsymbol{\mu}_{i} - \mathbf{r}_{i}\right\|_{2}^{2} + \|\text{diag}(\boldsymbol{\lambda}_{i})^{1/2}\boldsymbol{\mu}_{i}\|_{2}^{2}\right)$
            \label{step:bayesopinf-variance}
        \State $\boldsymbol{\Sigma}_{i} \gets \sigma_{i}^{*2}\left(\mathbf{D}\trp\mathbf{D} + \textup{diag}(\boldsymbol{\lambda}_{i})\right)^{-1}$
            \Comment{Covariance of operator row $\hat{\mathbf{o}}_{i}$.}
    \EndFor

    \LineComment{Online phase:~sample from operator posterior $p(\hat{\mathbf{O}}|\:\mathbf{D},\mathbf{R})$ and obtain ROM solutions.}
    \For{$\ell = 1, \ldots, N$}
        \For{$i = 1, \ldots, r$}
            \State $\hat{\mathbf{o}}_{i}^{(\ell)} \gets~\textrm{sample from}~\mathcal{N}(\boldsymbol{\mu}_{i},\boldsymbol{\Sigma}_{i})$
        \EndFor
        \State $\tilde{\mathbf{Q}}^{(\ell)} \gets~\textrm{integrate \eqref{eq:reducedorder} with}~\hat{\mathbf{O}}=[~\hat{\mathbf{A}}~\hat{\mathbf{H}}~\hat{\mathbf{B}}~\hat{\mathbf{c}}~] = [~\hat{\mathbf{o}}_{1}^{(\ell)}~\cdots~\hat{\mathbf{o}}_{r}^{(\ell)}~]\trp$ from $\hat{\mathbf{q}}_{0} = \hat{\mathbf{Q}}_{:,0}$ over $[t_0,t_f]$
    \EndFor
    \State \textbf{return}~$\tilde{\mathbf{Q}}^{(1)},\ldots,\tilde{\mathbf{Q}}^{(N)}$
        \Comment{ROM solutions corresponding to the operator samples.}
\EndProcedure
\end{algorithmic}
\caption{Bayesian Operator Inference (BayesOpInf).}
\label{alg:BayesOpInf}
\end{algorithm}

\vspace{3mm}
\noindent\textbf{Remark 4:} As the posterior $p(\hat{\vb{o}}_i |~ \vb{D},\vb{r}_i)$ follows a multivariate normal distribution, the noise-free posterior prediction of $\frac{\dd}{\dd t} \hat{q}_i(t)$ is a finite-rank Gaussian process, where $\hat{q}_{i}(t)$ is the $i$th entry of $\hat{\vb{q}}(t)$. Defined by the ROM form~\eqref{eq:reducedorder}, the posterior mean and covariance functions of $\frac{\dd}{\dd t} \hat{q}_i(t)$ take the form of
\begin{equation}
\begin{split}
    \mathbb{E}\left[\frac{\dd}{\dd t} \hat{q}_i(t) \middle\vert~\vb{D},\vb{r}_i\right]
    & = \vb{d}(\hat{\vb{q}}(t), \vb{u}(t))\trp\vb*{\mu}_i\,,
    \quad \text{and}\\
    \mathbb{C}ov\left[\frac{\dd}{\dd t} \hat{q}_i(t),\frac{\dd}{\dd t}  \hat{q}_i(t') \middle\vert~\vb{D},\vb{r}_i\right]
    & = \vb{d}(\hat{\vb{q}}(t), \vb{u}(t))\trp\vb*{\Sigma}_i \vb{d}(\hat{\vb{q}}(t'), \vb{u}(t'))\,,
\end{split}
\end{equation}
respectively,
where $\vb{d}(\hat{\vb{q}},\vb{u}) := [~\hat{\vb{q}}\trp~~(\hat{\vb{q}}\:\otimes\:\hat{\vb{q}})\trp~~\vb{u}\trp~~1~]\trp \in \mathbb{R}^{d(r,m)}$.
The same result can be obtained using Gaussian process regression \cite{rasmussen2006gaussian} with the following noise-corrupted prior:
\begin{equation}
\frac{\dd}{\dd t}  \hat{q}_i(t) \sim \mathcal{GP}\bigg(\vb{d}(\hat{\vb{q}}(t), \vb{u}(t))\trp\vb*{\beta}_i, \sigma_i^2\vb{d}(\hat{\vb{q}}(t), \vb{u}(t))\trp\mathrm{diag}(\vb*{\lambda}_i)^{-1}\vb{d}(\hat{\vb{q}}(t'), \vb{u}(t'))+\sigma_i^2\delta_\text{Dirac}(|t-t'|)\bigg)\,,
\end{equation}
in which $\delta_\text{Dirac}(\cdot)$ denotes the Dirac-delta function. This prior Gaussian process is consistent with the likelihood and prior definitions adopted in \eqref{eq:likelihood} and \eqref{eq:prior}, respectively, for the Bayesian inference \eqref{eq:post}.

\subsection{Equivalence to ridge regression}
\label{subsection:ridge}

Recall the posterior mean vector of $\hat{\vb{o}}_i$ given in \eqref{eq:musigma}. The correction term $\delta \vb*{\mu}_i= \vb*{\mu}_i - \vb*{\beta}_i$ can be expressed alternatively as
\begin{equation}\label{eq:ridge}
    \delta \vb*{\mu}_i
    = \underset{\vb*{\eta}\in \mathbb{R}^{d(r,m)}}{\textrm{arg\,min}}\left\{
        \|\vb{r}_i -\vb{D}(\vb*{\beta}_i +\vb*{\eta})\|_2^2
      + \|\mathrm{diag}(\vb*{\lambda}_i)^{\frac{1}{2}}\vb*{\eta}\|_2^2
    \right\}\,,
\end{equation}
implying that the posterior mean is the solution of a \emph{ridge regression} problem, i.e., a linear least-squares problem with a Tikhonov regularization (in this case $\|\mathrm{diag}(\vb*{\lambda}_i)^{\frac{1}{2}}\vb*{\eta}\|_2^2$). The use of Tikhonov regularization in OpInf has been reported in \cite{MHW2021regOpInfCombustion, SKHW2020romCombustion}. Thus another source of modeling uncertainty is introduced to the ROM by the Tikhonov regularization defined with OpInf. Based on the equivalence between the posterior mean of the inferred operators and ridge regression, it can be claimed that the uncertainty introduced by Tikhonov regularization is quantified through the proposed Bayesian framework as well.

On the other hand, the selection of the regularization parameters $\boldsymbol{\lambda}_{1},\ldots,\boldsymbol{\lambda}_{r}$ has been shown to be important in addressing ill-conditioning, over-fitting, or model inadequacy in OpInf \cite{MHW2021regOpInfCombustion,SKHW2020romCombustion}. The proposed Bayesian method in this work not only formulates a framework of uncertainty quantification, but also provides a new perspective for understanding the regularization selection process for deterministic OpInf. We address this further in Section \ref{sec:regularization}.

\section{Regularization selection}
\label{sec:regularization}

Without regularization, the OpInf point estimate~\eqref{eq:pointest} is represented with the Moore–Penrose inverse of matrix $\vb{D}$ and requires that $\vb{D}$ has full column rank, i.e., the Gram matrix $\vb{D}\trp\vb{D}$ is invertible. In practice, however, the condition number of $\vb{D}\trp\vb{D}$ may be large, which can lead to a poor estimate of the reduced operators and compromise the stability of the corresponding ROM. In this case, making use of a proper regularization can be critical for stabilizing the OpInf scheme and improving the robustness of the ROM \cite{MHW2021regOpInfCombustion,SKHW2020romCombustion}.

As specified in Section~\ref{subsection:ridge}, the posterior mean of inferred operators is equivalent to an estimator through ridge regression. The penalty coefficients $\vb*{\lambda}_i$ in the Tikhonov regularization correspond to the prior variances of the operators $\hat{\vb{o}}_i$. As given in \eqref{eq:maxlogmlhd}, the values of the hyperparameters $\vb*{\lambda}_i$ can be determined by an empirical Bayes method of maximizing marginal likelihood,
\begin{equation}\label{eq:marginal-likelihood}
\begin{split}
\mathcal{L}_i(\sigma_i^2,\vb*{\lambda}_i)  = &~ \log~\mathcal{N}\left(\vb{r}_i | ~\vb{D}\vb*{\beta}_i, \sigma_i^2\left(\vb{D}\mathrm{diag}(\vb*{\lambda}_i)^{-1}\vb{D}\trp+\vb{I}_k\right)\right) \\
 = &  -\frac{1}{2\sigma_i^2}\|\vb{r}_i-\vb{D}\vb*{\mu}_i \|_2^2-\frac{1}{2\sigma_i^2} \left\|\mathrm{diag}\left(\vb*{\lambda}_i\right)^{1/2}\delta\vb*{\mu}_i\right\|_2^2\\
 &  -\frac{1}{2}\log~\left|\mathrm{diag}(\vb*{\lambda}_i)+\vb{D}\trp\vb{D}\right| +\frac{1}{2}\log(\vb*{\lambda}_i)\trp\vb{1}_{d(r,m)} -\frac{k}{2}\log(\sigma_i^2)-\frac{k}{2}\log(2\pi)\,.
\end{split}
\end{equation}
In this sense, the empirical method provides guidance for determining the penalty parameters $\vb*{\lambda}_i$ in regularized OpInf.

Taking partial derivatives of $\mathcal{L}_i$ with respect to $\sigma_i^2$ and all the entries of $\vb*{\lambda}_i = [\lambda_{i,1},~\ldots,~\lambda_{i,d(r,m)}]\trp$, the extremum conditions are
\begin{equation*}
    \frac{\partial\mathcal{L}_i}{\partial \sigma_i^2}
    = 0
    \quad\Longrightarrow \quad
    \sigma_i^{*2}
    = \frac{\|\vb{r}_i-\vb{D}\vb*{\mu}_i\|_2^2+\left\|\mathrm{diag}\left(\vb*{\lambda}_i^*\right)^{1/2}\delta\vb*{\mu}_i\right\|_2^2}{k}\,,
\end{equation*}
which has already been given in \eqref{eq:noise}, and
\begin{equation}\label{eq:lambda}
    \frac{\partial\mathcal{L}_i}{\partial \lambda_{i,j}}
    = 0
    \quad \Longrightarrow\quad \lambda_{i,j}^{*}\left(\delta\mu_{i,j}^2+\Sigma_{i,jj}\right) = \sigma_i^{*2}\,,\quad 1\leq j \leq d(r,m)\,,
\end{equation}
in which $\delta\mu_{i,j}$ and $\Sigma_{i,jj}$ are the $j$-th entry of $\delta\vb*{\mu}_i$ and the $j$-th diagonal entry of $\vb*{\Sigma}_i$, respectively.\footnote{The evaluation of $\frac{\partial\mathcal{L}_i}{\partial \lambda_{i,j}}$ uses the identity $\partial_\theta (\log\det \vb{A}(\theta))=\tr(\vb{A}(\theta)^{-1}\partial_\theta\vb{A}(\theta))$, in which $\vb{A}$ is an invertible matrix parametrized by $\theta$.}
Note that \eqref{eq:lambda} is an implicit representation of $\vb*{\lambda}_i^{*}$, the penalty coefficients of Tikhonov regularization.

For each reduced system mode $i=1,\ldots,r$, the problem of maximum marginal likelihood in \eqref{eq:marginal-likelihood} is non-convex with $1 + d(r,m)$ unknowns. Solving such a problem using gradient-based methods is possible, but inconvenient for moderately sized $r$ and not guaranteed to find a global optimum. To alleviate complexity, we impose the constraint $\lambda_{i,j} = \lambda_{i}$ for all $1\leq j \leq d(r,m)$, hence $\boldsymbol{\lambda}_{i} = [~\lambda_{i}~\lambda_{i}~\cdots~\lambda_{i}~]\trp\in\mathbb{R}^{d(r,m)}$.
In this case, the regularization term in \eqref{eq:leastsqr}, as well as in steps~\ref{step:bayesopinf-mean} and~\ref{step:bayesopinf-variance} of Algorithm~\ref{alg:BayesOpInf}, is an $L_2$-regularization with penalty coefficient $\lambda_{i} > 0$, i.e.,
\begin{align}
    \mathcal{P}_{i}(\boldsymbol{\mu}_{i})=\mathcal{P}_{i}(\vb*{\beta}_i+\delta\boldsymbol{\mu}_{i})
    = \left\|\textrm{diag}(\boldsymbol{\lambda}_{i})^{1/2}\delta\boldsymbol{\mu}_{i}\right\|_{2}^{2}
    = \lambda_{i}\left\|\delta\boldsymbol{\mu}_{i}\right\|_{2}^{2}.
\end{align}
This means that the operator entries corresponding to the reduced system mode $\hat{q}_{i}$, i.e., the $i$th rows of $\hat{\mathbf{c}}$, $\hat{\mathbf{A}}$, $\hat{\mathbf{H}}$, and $\hat{\mathbf{B}}$, are regularized equally, while operator entries corresponding to different system modes are regularized separately.
Summing \eqref{eq:lambda} over all $1\leq j \leq d(r,m)$ yields
\begin{equation}
    \lambda_i^*(\|\delta\vb*{\mu}_i\|_2^2+\tr(\vb{\Sigma}_i))
    = d(r,m) \sigma_i^{*2}\,,
\end{equation}
in which $\tr(\vb{\Sigma}_i)$ can alternatively be written as follows \cite{rasmussen2006gaussian}:
\begin{equation}
    \tr(\vb*{\Sigma}_i)
    = \sigma_i^{*2}\sum_{l=1}^{d(r,m)}\frac{1}{\lambda_i^*+g_l}\,,
\end{equation}
where $g_1,\ldots,g_{d(r,m)}\ge 0$ are the non-negative eigenvalues of the Gram matrix $\vb{D}\trp\vb{D}$. This representation of $\tr(\vb{\Sigma}_i)$ has the computational advantage of avoiding explicit matrix inversion in \eqref{eq:musigma} to evaluate $\boldsymbol{\Sigma}_{i}$. We then have
\begin{equation}
    \label{eq:lambda2}
    \lambda^*_i
    = \frac{\sigma_i^{*2}}{\|\delta\vb*{\mu}_i\|_2^2}\underbrace{\sum_{l=1}^{d(r,m)}\frac{g_l}{\lambda_i^{*}+g_l} }_{:=\gamma_i}
    = \frac{\gamma_i\sigma_i^{*2}}{\|\delta\vb*{\mu}_i\|_2^2}\,,
\end{equation}
in which $0< \gamma_i < d(r,m)$. The following equality shows the connections among the noise $\sigma_i^{*2}$ and the two terms in the loss function of the ridge regression \eqref{eq:ridge}:
\begin{equation}
    \frac{\lambda_i^*\|\delta\vb*{\mu}_i\|_2^2}{\gamma_i}
    = \frac{\|\vb{r}_i-\vb{D}\vb*{\mu}_i\|_2^2}{k-\gamma_i}
    = \sigma_i^{*2}\,.
\end{equation}

Substituting \eqref{eq:noise} into \eqref{eq:lambda2}, ${\lambda}_i^*$ is expressed by the posterior mean $\delta\vb*{\mu}_i$ and covariance matrix $\vb*{\Sigma}_i$, both being dependent on ${\lambda}_i^*$ as well. Thus an implicit representation of ${\lambda}_i^*$ is given by ${\lambda}_i^*=\mathcal{F}({\lambda}_i^*)$, where $\mathcal{F}$ denotes the right-hand side of \eqref{eq:lambda2}. In this work, we suggest a fixed-point iterative strategy, $({\lambda}_i^*)_{s+1}=\mathcal{F}(({\lambda}_i^*)_s)$, $s = 0, 1, \ldots$, as a practical approach to simultaneously satisfy \eqref{eq:noise} and \eqref{eq:lambda2} and approximately maximize \eqref{eq:marginal-likelihood}. The initial values $\{({\lambda}_i^*)_{0}\}_{i=1}^{r}$ can be set to a selection of regularization coefficients that result in a stable ROM, but which are not necessarily optimal. This iterative regularization update procedure is given by Algorithm~\ref{alg:BayesOpInf-Case2}, which considers a zero prior mean of reduced operators, i.e., $\vb*{\beta}_i = \vb{0}$ and $\vb*{\mu}_i = \delta\vb*{\mu}_i$. We note that while Algorithm~\ref{alg:BayesOpInf-Case2} is computationally efficient, it does not include stability information for the reduced-order ODE system~\eqref{eq:reducedorder} because of its completely non-intrusive nature.

\begin{algorithm}
\begin{algorithmic}[1]
\Procedure{BayesOpInfRegularization}{
    Data matrix $\mathbf{D}\in\mathbb{R}^{k\times d(r,m)}$,
    projected time derivatives $\mathbf{R} = [~\mathbf{r}_{1}~\cdots~\mathbf{r}_{r}~]\trp\in\mathbb{R}^{r\times k}$,
    initial regularization parameters $\lambda_{1},\ldots,\lambda_{r} \ge 0$,
    convergence tolerance $\varepsilon > 0$
    }

    \State $g_{1},\ldots,g_{d(r,m)} \gets \texttt{eig}(\mathbf{D}\trp\mathbf{D})$
        \Comment{Get eigenvalues of the data Gram matrix.}
    \State $\boldsymbol{\lambda}^{(0)} \gets [~\lambda_1~\cdots~\lambda_{r}~]$
    \For{$\ell = 0, 1, \ldots$}
        \For{$i= 1, \ldots, r$}
            \State $\boldsymbol{\mu}_{i} \gets \arg\min_{\boldsymbol{\mu}_{i}}\left\|\mathbf{D}\boldsymbol{\mu}_{i} - \mathbf{r}_{i}\right\|_{2}^{2} + \lambda_{i}\left\|\boldsymbol{\mu}_{i}\right\|_{2}^{2}$
                \Comment{Update operator mean.}
            \State $\sigma_{i}^{*2} \gets \frac{1}{k}\left(\left\|\mathbf{D}\boldsymbol{\mu}_{i} - \mathbf{r}_{i}\right\|_{2}^{2} + \lambda_{i}\left\|\boldsymbol{\mu}_{i}\right\|_{2}^{2}\right)$
                \Comment{Update operator variance.}
            \State $\lambda_{i}^{*} \gets \displaystyle\frac{\sigma_{i}^{*2}}{\|\boldsymbol{\mu}_{i}\|_{2}^{2}}\displaystyle\sum_{l=1}^{d(r,m)}\frac{g_{l}}{\lambda_{i} + g_{l}}$
                \Comment{Update regularization.}
        \EndFor
        \State $\boldsymbol{\lambda}^{(\ell+1)} \gets [~\lambda_1^{*}~\cdots~\lambda_{r}^{*}~]$
        \If{$\|\boldsymbol{\lambda}^{(\ell+1)} - \boldsymbol{\lambda}^{(\ell)}\|_2 \big/ \|\boldsymbol{\lambda}^{(\ell)}\|_2 < \varepsilon$}
            \Comment{Check for convergence.}
            \label{step:relativechange-convergence}
            \State \textbf{return}~$\lambda_{1}^{*},\ldots,\lambda_{r}^{*}$
        \EndIf
        \State $\lambda_{1},\ldots,\lambda_{r} \gets \lambda_{1}^{*},\ldots,\lambda_{r}^{*}$
        \Comment{Continue iterating if not converged.}
    \EndFor
\EndProcedure
\end{algorithmic}
\caption{Iterative regularization update for BayesOpInf with $\lambda_{i,j} = \lambda_{i}$.}
\label{alg:BayesOpInf-Case2}
\end{algorithm}

\section{Numerical examples}
\label{sec:results}

Section \ref{sec:euler} demonstrates the proposed method on one-dimensional Euler equations where the state solution data are polluted by Gaussian noise.
In Section \ref{sec:combustion}, we examine a combustion application in which the regularizing terms play a crucial role for learning stable ROMs from data.
The code for these experiments can be found at \url{https://github.com/Willcox-Research-Group/ROM-OpInf-Combustion-2D/tree/cmame2022}.

\subsection{Noised Euler Equations}
\label{sec:euler}

Consider the conservative one-dimensional compressible Euler equations for an ideal gas,
\begin{align}
    \label{eq:euler-conservative}
    \frac{\partial}{\partial t}\left[\rho\right]
    &= -\frac{\partial}{\partial x}\left[\rho u\right],
    &
    \frac{\partial}{\partial t}\left[\rho u\right]
    &= -\frac{\partial}{\partial x}\left[\rho u^2 + p\right],
    &
    \frac{\partial}{\partial t}\left[\rho e\right]
    &= -\frac{\partial}{\partial x}\left[(\rho e + p)u\right],
\end{align}
where $u$ is velocity [m/s], $\rho$ is density [kg/m$^3$], $p$ is pressure [Pa], $\rho u$ is specific momentum [kg/m$^2$s], and $\rho e$ is total energy [kgJ/m$^3$].
The state variables are related via the ideal gas law $\rho e = \frac{p}{\gamma - 1} + \frac{1}{2}\rho u^2$, where $\gamma = 1.4$ is the heat capacity ratio.
Our goal is to learn a reduced model of \eqref{eq:euler-conservative} from noisy observations of the conservative variables $(\rho, \rho u, \rho e)$ and to quantify the uncertainties induced by the noise. The ideal gas law enables a change of variables $\tau:(\rho, \rho u, \rho e) \mapsto (u, p, 1/\rho)$ that transforms \eqref{eq:euler-conservative} into the purely quadratic system
\begin{align}
    \label{eq:euler-lifted}
    \frac{\partial u}{\partial t}
    &= -u \frac{\partial u}{\partial x} - \zeta\frac{\partial p}{\partial x},
    &
    \frac{\partial p}{\partial t}
    &= -\gamma p \frac{\partial u}{\partial x} - u\frac{\partial p}{\partial x},
    &
    \frac{\partial \zeta}{\partial t}
    &= -u \frac{\partial\zeta}{\partial x} + \zeta\frac{\partial u}{\partial x},
\end{align}
where $\zeta = 1/\rho$ is the specific volume [m$^3$/kg].
In the model~\eqref{eq:euler-lifted}, every term is the product of exactly two state variables $(u,p,\zeta)$ or their spatial derivatives, hence for this problem we seek a reduced model with quadratic terms only:
\begin{align}
    \label{eq:rom-quadratic}
    \frac{\textup{d}}{\textup{d}t}\hat{\mathbf{q}}(t)
    = \hat{\mathbf{H}}[\hat{\mathbf{q}}(t)\otimes\hat{\mathbf{q}}(t)].
\end{align}
The OpInf operator matrix of \eqref{eq:matrix-definitions} in this case is $\hat{\mathbf{O}} = \hat{\mathbf{H}} \in\mathbb{R}^{r\times d(r)}$, with column dimension $d(r) = r(r+1)/2$.

\begin{figure}
    \centering
    \includegraphics[width=\textwidth]{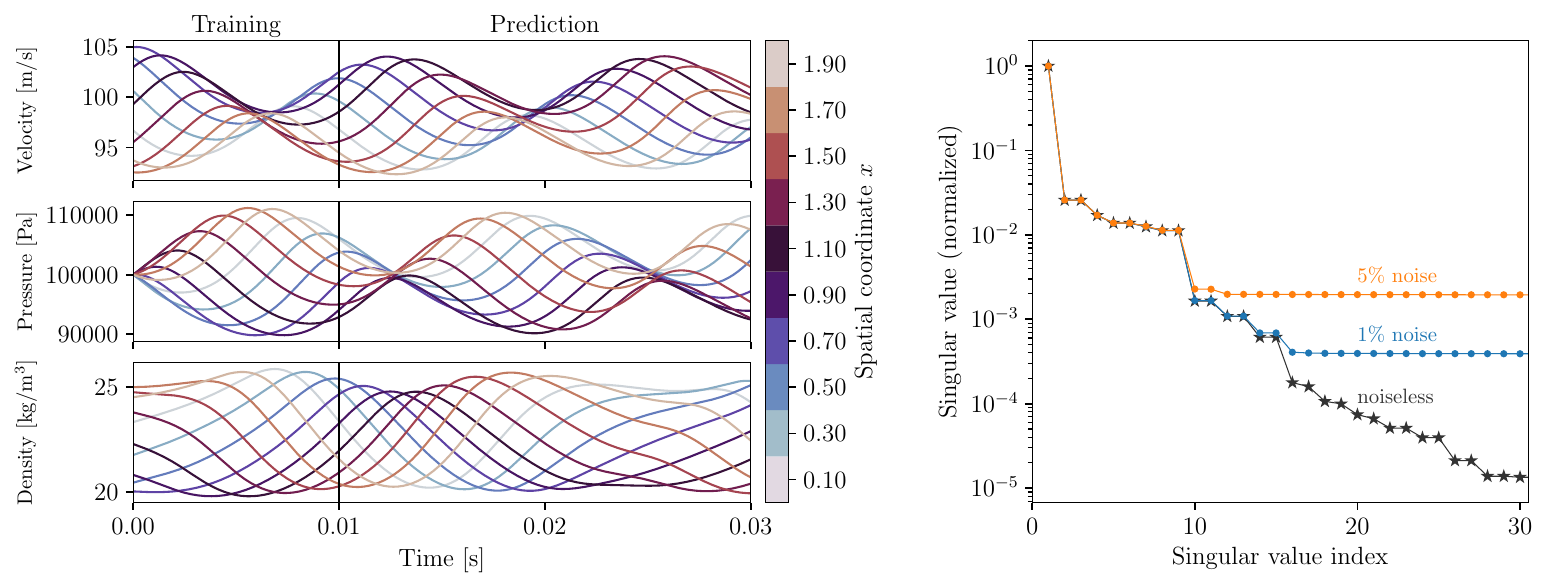}
    \vspace{-.75cm}
    \caption{One of the sixty-four training trajectories at various spatial coordinates before noise corruption (left) and singular value decay of the entire training set with various noise levels (right). Training data consists of noised snapshots for $0 \le t < 0.01$.}
    \label{fig:euler-training}
\end{figure}

We consider the spatial domain $[0,2)$ and the time domain $[t_0,t_f] = [0, 0.03]$, with periodic boundary conditions imposed at the spatial boundaries.
Training data is generated over the shorter time domain $[0, 0.01]$ for multiple initial conditions: the initial pressure is set to $10^5$~Pa everywhere, and initial velocity and density profiles are constructed via cubic spline interpolation by fixing values at the spatial points $x \in \{0, \frac{2}{3}, \frac{4}{3}\}$ to either $95$ or $105$~m/s for velocity and either $20$ or $24$~kg/m$^3$ for density. We therefore have $2^{3\times 2} = 64$ initial conditions (three interpolation nodes with two choices of values for two variables).
For each initial condition, we solve the conservative system~\eqref{eq:euler-conservative} by discretizing the spatial domain via first-order finite differences with spatial resolution $\delta x = 0.01$, then stepping forward in time via the first-order explicit Euler iteration with time step $\delta t = 10^{-5}$. The result is $k = 64\times 1000 = 64{,}000$ total training snapshots, each of size $n = 200\times 3 = 600$. See Figure~\ref{fig:euler-training} for an example trajectory in the training set.
Before processing the snapshot data for learning, each snapshot entry is polluted with Gaussian noise proportional to the range of the associated variable to represent observation error. For example, to each snapshot entry corresponding to density we add error drawn from $\mathcal{N}(0, \varsigma^{2})$ where $\varsigma = \xi(\max_{i,j}\{\rho_{i,j}\} - \min_{i,j}\{\rho_{i,j}\})$,  $\xi = 5\%$ is the noise level, and $\rho_{i,j}$ is the $i$th entry of the $j$th density training snapshot.
The noisy snapshots are then transformed from the conservative states of \eqref{eq:euler-conservative} to the states of the quadratic system \eqref{eq:euler-lifted} via the mapping $\tau$ and non-dimensionalized, yielding training data $\mathbf{q}_{0},\ldots,\mathbf{q}_{k-1}\in\mathbb{R}^{n}$ to be used in OpInf. Note that the noise in the lifted variables is not necessarily Gaussian since $\tau$ is a nonlinear mapping.

\begin{figure}
    \centering
    \includegraphics[width=\textwidth]{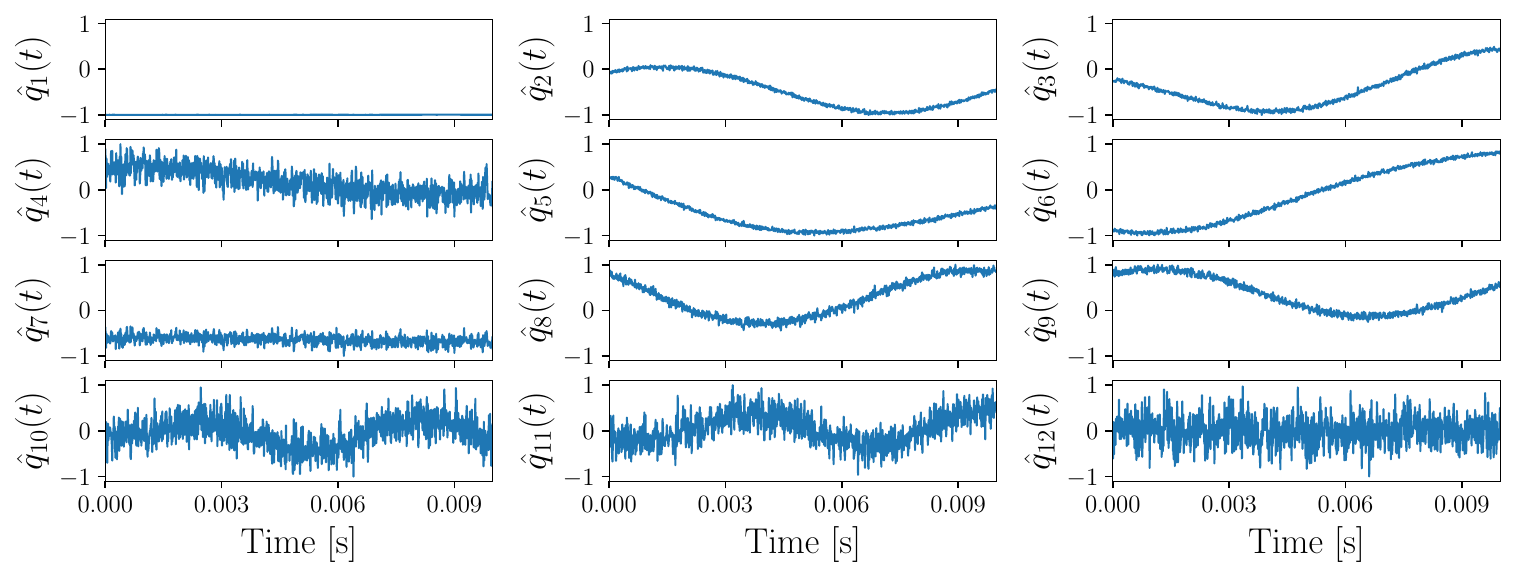}
    \vspace{-.75cm}
    \caption{The first several modes of the (noisy) projected training data, normalized to the range $[-1,1]$, for the training trajectory shown in Figure~\ref{fig:euler-training}. The POD projection successfully filters much, but not all, of the noise on the first $r = 9$ modes.}
    \label{fig:euler-training-projected}
\end{figure}

\begin{figure}[!ht]
    \centering
    \includegraphics[width=\textwidth]{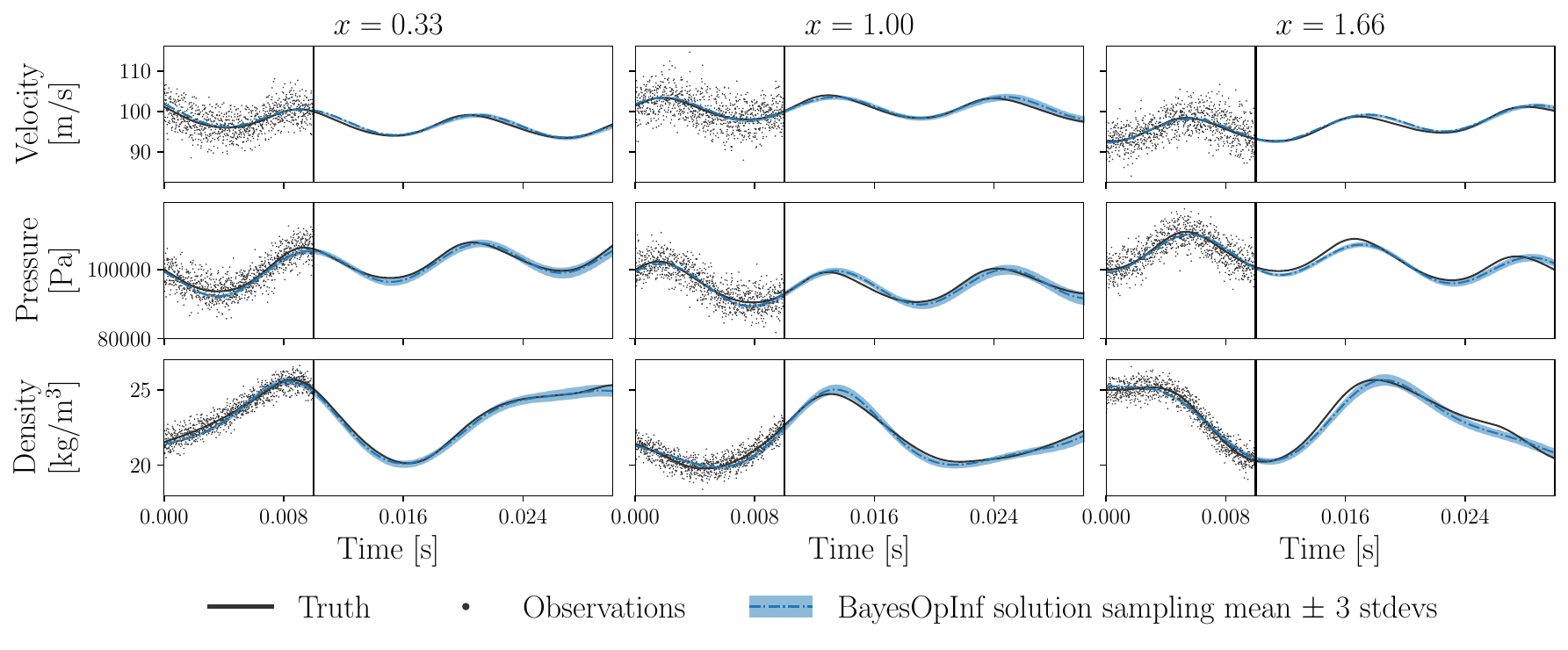}
    \vspace{-.75cm}
    \caption{Traces of full-order and ROM solutions to the Euler problem at three spatial locations, with the same initial condition as in Figure~\ref{fig:euler-training}. The noisy training data (noise level $\xi = 5\%$) is also shown as individual points over the training interval. The BayesOpInf ROM uses $r = 9$ modes, with the regularization chosen via Algorithm \ref{alg:BayesOpInf-Case2}, and 100 draws are made from the posterior distribution to sample the mean and standard deviations.
    }
    \label{fig:euler-traces}
\end{figure}

From the lifted, noisy snapshots we compute a POD basis $\mathbf{V}\in\mathbb{R}^{n\times r}$ and project the training data to the corresponding low-dimensional linear subspace as $\hat{\mathbf{q}}_{j} = \mathbf{V}\trp\mathbf{q}_{j}$. The dimension $r = 9$ is selected based on where the decay of the singular values flattens out due to the noise pollution as seen in Figure~\ref{fig:euler-training}. It has been widely observed that POD has a smoothing effect (see, e.g., \cite{del2008proper, epps2010PODthresh, venturi2006PODperturb}). The projected training data, an example of which is shown in Figure~\ref{fig:euler-training-projected}, exhibits relatively little noise when compared to the observed training data (see also Figure~\ref{fig:euler-traces}) up to $r = 9$ modes. However, the OpInf regression requires an estimate of the time derivatives $\dot{\mathbf{q}}_{j}$, which will be highly inaccurate if computed as finite differences of the noisy projected data. We therefore apply the noise-resistant strategy of \cite{DeBrabanter2013localddt}, which uses implicit local polynomial regression to estimate derivatives with minimal variance, to estimate each $\dot{\mathbf{q}}_{j}$. Therefore, uncertainty in the OpInf regression stems from \emph{i}) noise in the observed data, \emph{ii}) the truncated POD representation (model form error), and \emph{iii}) error in the estimates of the time derivatives of the observed data.

The noise in the training data necessitates some regularization in order to achieve stable reduced models, though in this example among stable models the learned OpInf ROM is not highly sensitive to the choice of regularization. We initially set $\lambda_{i,j} = 50$ for all $i,j$ and use Algorithm~\ref{alg:BayesOpInf-Case2} to iterate until the relative change (as measured in step~\ref{step:relativechange-convergence}) drops below $\varepsilon = 0.1\%$, which occurs after only nine iterations.
We then draw $100$ samples of $\hat{\mathbf{O}}$ from the posterior operator distribution defined in \eqref{eq:final-operator-posterior} and integrate the corresponding ROMs over the full time domain $[t_0,t_f]$ for each of the (noisy) initial conditions in the training set. Figure~\ref{fig:euler-traces} shows the sample mean of these draws and three standard deviations from the sample mean at several spatial locations for the trajectory shown in Figure~\ref{fig:euler-training}; the results are representative for the set of all initial conditions considered. Note that the width of three standard deviations, which slightly increases with time, is relatively small compared to the observed data. This is because the uncertainty band shows a credible interval, which quantifies the accuracy of model fits, i.e., the learned posterior distribution describes the uncertainty relative to the true solution. Note that, over the training interval $[0,0.01]$, the credible interval is by definition not taking into account the measurement noise in the snapshot data\footnote{The band that takes both the propagated uncertainty from the posterior model and the measurement errors (noise) into account is often referred to as a prediction interval, see chapter 9 of \cite{smith2013uncertainty}.}.

\begin{figure}
    \centering
    \includegraphics[width=\textwidth]{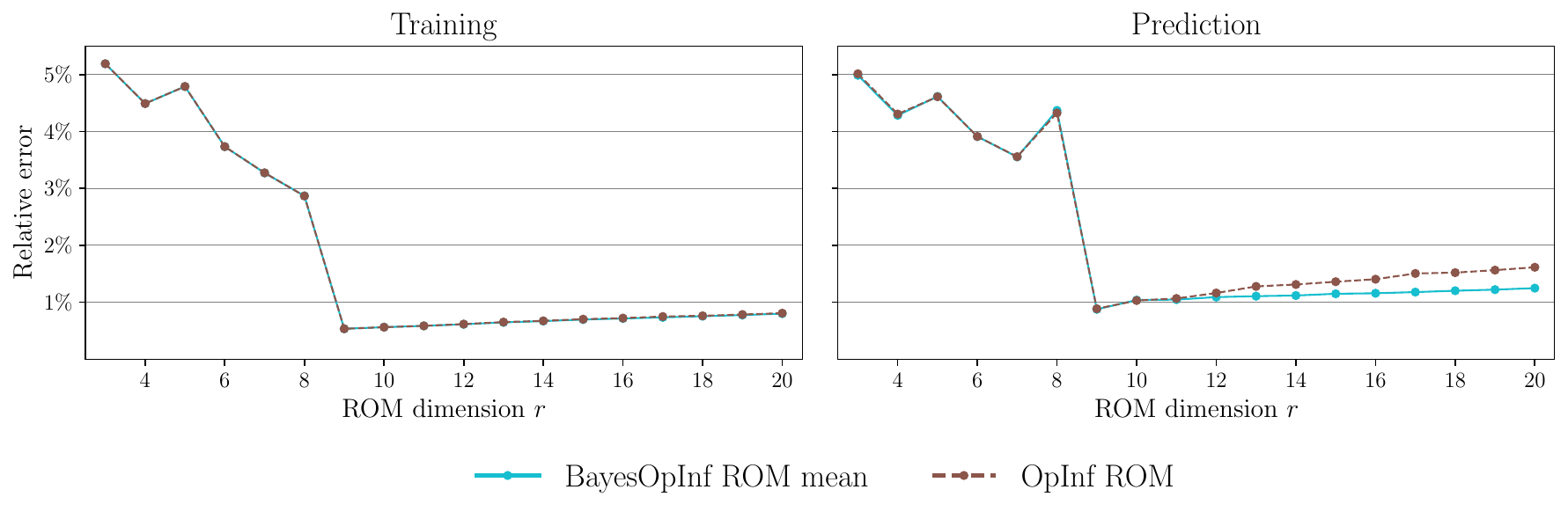}
    \vspace{-.75cm}
    \caption{Average relative error over the set of 64 initial conditions for increasing basis sizes, over the training regime (left) and the prediction regime (right), for ROMs where the regularization is chosen via Algorithm~\ref{alg:BayesOpInf-Case2} (BayesOpInf) and by Algorithm~\ref{alg:OpInf-Reg} (OpInf).
    }
    \label{fig:euler-basisVerror}
\end{figure}

The regularization selection strategy of Algorithm~\ref{alg:BayesOpInf-Case2} only approximates the maximizers of the marginal likelihood~\eqref{eq:marginal-likelihood}, however the cost of each iteration in Algorithm \ref{alg:BayesOpInf-Case2} is dominated only by the OpInf regression (determining each $\boldsymbol{\mu}_{i}$ and $\boldsymbol{\Sigma}_i$). To compare this strategy to existing methods, Figure~\ref{fig:euler-basisVerror} repeats the experiments described above with basis sizes $r = 3, 4, \ldots, 20$ using Algorithm \ref{alg:BayesOpInf-Case2} to select the regularization, and computes the average relative error of the ROM defined by the posterior operator mean $\boldsymbol{\mu}_{i}$ over the set of all initial conditions for both the training regime $t \in [0, 0.01)$ (with the error computed with respect to the true trajectories, not the noisy observations) and the prediction regime $t\in [0.01, 0.03]$. The same results are shown for the more involved regularization selection strategy of~\cite{MHW2021regOpInfCombustion} detailed in Algorithm~\ref{alg:OpInf-Reg} (see Appendix~2). This strategy chooses a scalar regularization value (i.e., $\lambda_{i,j} = \lambda$ for all $i,j$) by solving an optimization problem based on two criteria: \emph{i}) that the resulting ROM remains stable throughout the full time horizon $t \in [0,0.03]$ with each of the considered initial conditions, and \emph{ii}) that the reconstruction error of the training data is minimized. Note that Algorithm~\ref{alg:OpInf-Reg} chooses the scalar regularization that minimizes the difference between the integrated ROM state and the \emph{noisy} training data. Figure~\ref{fig:euler-basisVerror} shows that the errors for each selection method are highly similar, especially over the training regime, but that the BayesOpInf approach is slightly more accurate in the prediction regime for $r > 11$. This is likely because Algorithm \ref{alg:BayesOpInf-Case2} treats the regularization of each dynamic mode $\hat{q}_{i}$ separately, while Algorithm~\ref{alg:OpInf-Reg} uses a single global regularizer for all terms of the learned operators. The key takeaway is that when the stability of the learned model is not highly sensitive to the regularization hyperparameters, Algorithm \ref{alg:BayesOpInf-Case2} provides a lightweight option for selecting the regularization. The next section details an example where the choice of regularization is critical for model stability, in which case Algorithm~\ref{alg:BayesOpInf-Case2} is a less appropriate option for choosing the $\boldsymbol{\lambda}_{i}$.

\subsection{Single-injector Combustion Process}
\label{sec:combustion}

We now consider a combustion application on a two-dimensional domain in which model form error, and hence the choice of regularization, plays a significant role in the OpInf procedure. This problem presents a challenging application for model reduction and has been used as a test problem for both intrusive methods \cite{HDKM2018combustion-is-hard,HDM2019poddeim-robustness,HXDM2018rocketrom-poddeim} and non-intrusive OpInf \cite{jain2021performance,MHW2021regOpInfCombustion,SKHW2020romCombustion}.
The governing dynamics are the nonlinear conservation laws
\begin{align}
    \label{eq:combustion-conservative}
    \frac{\partial\vec{q}_\textrm{c}}{\partial t}
    + \nabla \cdot(\vec{K} - \vec{K}_{\textrm{v}})
    &= \vec{S},
\end{align}
where
$
    \vec{q}_\textrm{c}
    = (\rho, \rho v_x, \rho v_y, \rho e, \rho Y_1, \rho Y_2, \rho Y_3, \rho Y_4)
$
are the conservative variables, $\vec{K}$ is the inviscid flux, $\vec{K}_\textrm{v}$ is the viscous flux, and $\vec{S}$ contains source terms.
Here $\rho$ is the density $[\frac{\textrm{kg}}{\textrm{m}^3}]$,
$v_x$ and $v_y$ are the $x$ and $y$ velocities $[\frac{\textrm{m}}{\textrm{s}}]$,
$e$ is the total energy $[\frac{\textrm{J}}{\textrm{m}^3}]$,
and $Y_1$, $Y_{2}$, $Y_{3}$, $Y_{4}$ are the mass fractions for the chemical species CH$_4$, O$_2$, H$_2$O, and CO$_2$, respectively, which follow the global one-step reaction
$
    \textrm{CH}_4 + 2\textrm{O}_2
    \to
    \textrm{CO}_2 + 2\textrm{H}_2\textrm{O}
$
\cite{WD1981oxidation}. The top and bottom walls of the combustor are assigned no-slip boundary conditions, the upstream boundary is constrained to have a constant mass flow at the inlets, and at the downstream end a non-reflecting boundary condition is imposed to maintain the chamber pressure via
\begin{align}
    p_\textrm{back}(t)
    &= p_{\textrm{back,ref}}\left(1 + 0.1\sin(2\pi f t)\right),
    \label{eq:input-function}
\end{align}
where $p_{\textrm{back,ref}} = 10^6$~Pa and $f = 5{,}000$~Hz.
See \cite{HHSFAMT2015gems} for more details on the governing equations.

\begin{figure}
    \centering
    \includegraphics[width=.8\textwidth]{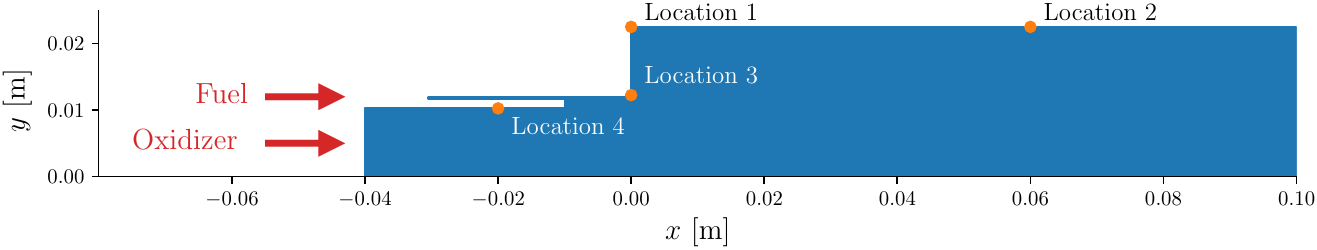}
    \vspace{-.25cm}
    \caption{Two-dimensional domain of the single-injector combustor. The solution is monitored pointwise at sensor locations $1$--$4$ and reported in Figures \ref{fig:combustion-pressure}--\ref{fig:combustion-timeVspread-xvelocity}.}
    \label{fig:combustion-domain}
\end{figure}

As observed in \cite{SKHW2020romCombustion}, the governing conservative equations~\eqref{eq:combustion-conservative} take on a nearly quadratic form when written in a specific volume form with variables
$
    \vec{q} = (p, v_x, v_y, T, \zeta, c_1, c_2, c_3, c_4),
$
where $p$ is the pressure $[\textrm{Pa}]$,
$T$ is the temperature $[\textrm{K}]$,
$\zeta = 1/\rho$ is the specific volume $[\frac{\textrm{m}^3}{\textrm{kg}}]$,
and $c_1,\ldots,c_4$ are the species molar concentrations $[\frac{\textrm{kmol}}{\textrm{m}^3}]$ given by $c_\ell = \rho Y_\ell / M_\ell$ with $M_\ell$ the molar mass of the $\ell$th species $[\frac{\textrm{g}}{\textrm{mol}}]$.
In \cite{MHW2021regOpInfCombustion}, it was shown that including pressure, specific volume, and temperature led to improved ROM performance, although only two of the three are needed to fully specify the high-fidelity model, with the equation of state defining the third. Hence, we construct a fully quadratic ROM of the form~\eqref{eq:reducedorder}, with the understanding that this structure only approximates the structure of the governing equations.

In this experiment, we consider a single set of initial conditions and the corresponding trajectory. Our objective is to predict in time beyond the available training data. We solve \eqref{eq:combustion-conservative} via the finite-volume based General Equation and Mesh Solver (GEMS), see \cite{HHSFAMT2015gems} for details. The spatial domain is discretized with $n_x = 38{,}523$ cells (hence each snapshot has $8n_x = 308{,}184$ entries) and the solution is computed for $50{,}000$ time steps with step size $\delta t = 10^{-7}$~s, from $t_0 = 0.015$~s to $t_{f} = 0.020$~s. We retain the first $k = 20{,}000$ snapshots for training and reserve the remaining $30{,}000$ for testing. The training snapshots are lifted to the specific volume variables $\vec{q}$, then scaled to $[-1,1]$ by first subtracting out the mean profile in the temperature, pressure, and specific volume, then normalizing each variable by its maximum absolute entry. These processed training snapshots $\mathbf{q}_{0},\ldots,\mathbf{q}_{k-1}$ have dimension $n = 9n_x = 346{,}000$. We compute the POD basis $\mathbf{V}$ of the training snapshots using a randomized singular-value decomposition \cite{HMPT2011rNLA} and retain $r = 38$ modes based on the associated singular value decay. As the training data is smooth (i.e., not corrupted by observational noise), we estimate time derivatives of the projected training data $\hat{\mathbf{q}}_{j} = \mathbf{V}\trp\mathbf{q}_{j}$ with fourth-order finite differences. This estimation introduces a further source of error, though the predominant error is due to the model misspecification by assuming a quadratic ROM form.

\begin{figure}
    \centering
    \includegraphics[width=\textwidth]{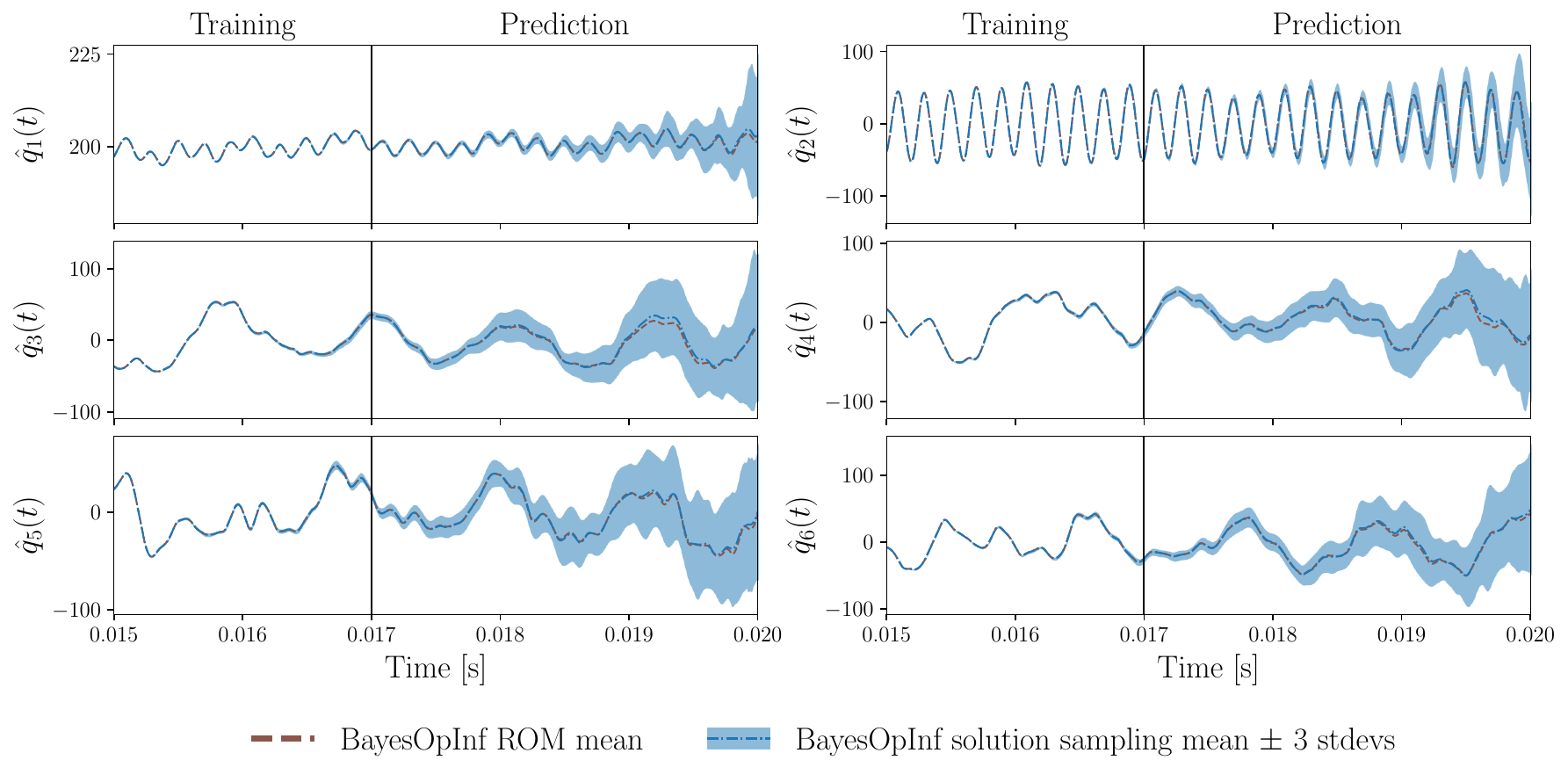}
    \vspace{-.75cm}
    \caption{Evaluations of the first six modes (i.e., the first six entries of $\hat{\mathbf{q}}(t)$) of $100$ draws from the posterior distribution of the Bayesian OpInf ROM. The dashed lines show the trajectory of the ROM with posterior mean operators ($\hat{\mathbf{o}}_{i} = \boldsymbol{\mu}_{i}$); the shaded areas are the regions within three standard deviations of the sample mean. Training data is observed for $t \in [0.015, 0.017)$.}
    \label{fig:combustion-modes}
\end{figure}

Previous applications of OpInf to this problem have shown that regularization is essential to obtain stable ROMs \cite{MHW2021regOpInfCombustion,SKHW2020romCombustion}. We therefore utilize the regularization selection strategy of \cite{MHW2021regOpInfCombustion} (Algorithm~\ref{alg:OpInf-Reg}) to choose the regularization parameters $\lambda_{i,j}$. In this context, we select two scalars: one to penalize the terms of the non-quadratic operators $\hat{\mathbf{c}}$, $\hat{\mathbf{A}}$, and $\hat{\mathbf{B}}$, and one for $\hat{\mathbf{H}}$ alone (see Appendix~2 for details). Thus the posterior mean $\boldsymbol{\mu}_{i}$ of the operator $\hat{\mathbf{o}}_{i}$ given in \eqref{eq:final-operator-posterior} is equivalent to the ROM inference studied in \cite{MHW2021regOpInfCombustion}. We additionally draw $100$ operator samples from \eqref{eq:final-operator-posterior} and integrate the resulting ROMs to quantify the uncertainties in the learned model. The computational cost of integrating $100$ ROMs is trivial compared to the high-dimensional code: computing the $50{,}000$ unprocessed snapshots with GEMS costs approximately $1{,}000$ CPU hours, while constructing and integrating a single ROM from the posterior takes $\sim 0.3$ CPU seconds. Even when the mean operators $\boldsymbol{\mu}_{i}$ of the posterior distribution~\eqref{eq:final-operator-posterior} define a stable ROM, there is no guarantee that operators drawn from the distribution will result in stable ROMs. In this experiment, $96$ of the $100$ posterior draws yield stable ROMs. This is an indicator of the complexity of the problem at hand and highlights the constraint-free nature of the posterior formulation.

Figure~\ref{fig:combustion-modes} plots the dominant modes of the reduced-order state as a function of time together with the deviation defined by the posterior draws. Both the trajectory corresponding to the mean operators $\boldsymbol{\mu}_{i}$ and the sample mean of the $96$ stable draws are shown, but they differ very little. We emphasize that, while the entries of the model operators are assumed to follow a Gaussian distribution, there is no reason for the integrated trajectories of the corresponding quadratic model to be Gaussian in nature. Notably, the uncertainty bands remain tight until the end of the training regime, and widen increasingly with time. The band is thinnest in the second mode, which corresponds to the pressure signal imparted by \eqref{eq:input-function}.

\begin{figure}[htbp]
    \centering
    \includegraphics[width=\textwidth]{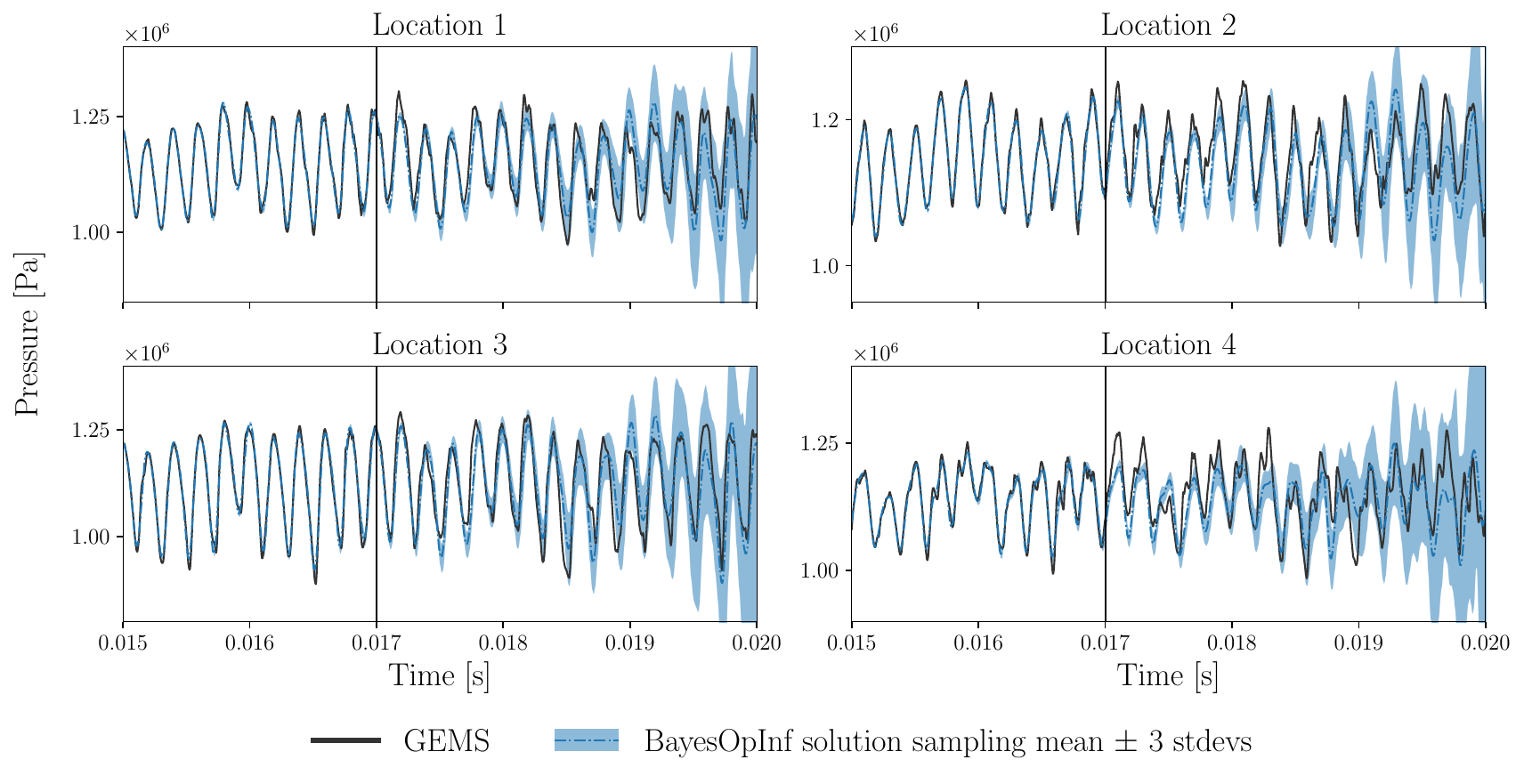}
    \vspace{-.75cm}
    \caption{Full-order GEMS data and Bayesian OpInf ROM reconstructions of the pressure at the four sensor locations marked in Figure~\ref{fig:combustion-domain}. The shaded areas are the regions within three standard deviations of the sample mean, using the same $100$ draws as in Figure~\ref{fig:combustion-modes}.
    Training data is observed for $t \in [0.015, 0.017)$.}
    \label{fig:combustion-pressure}
    \bigskip\bigskip
    \includegraphics[width=\textwidth]{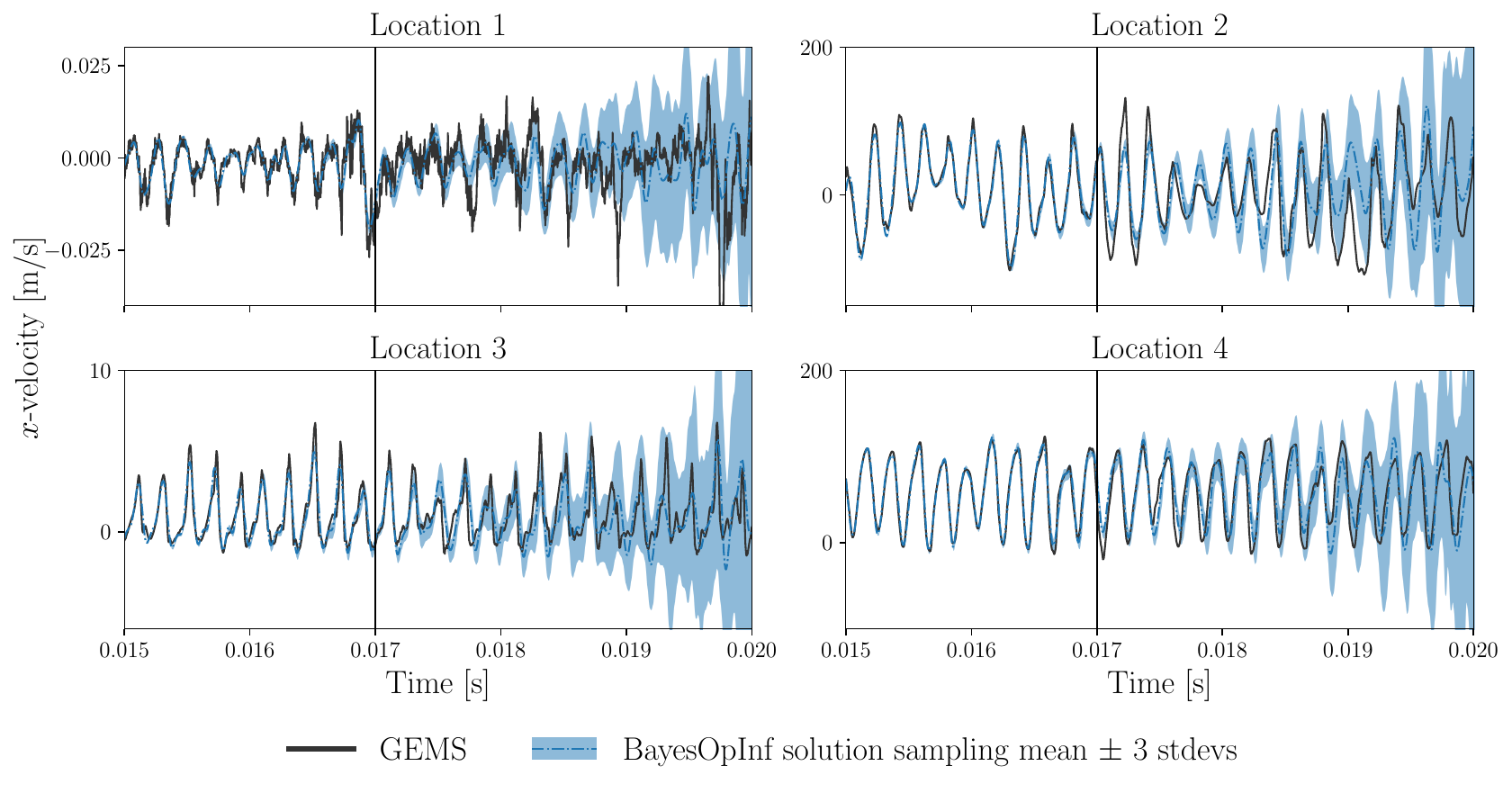}
    \vspace{-.75cm}
    \caption{Full-order GEMS data and Bayesian OpInf ROM reconstructions of the horizontal velocity at the four sensor locations marked in Figure~\ref{fig:combustion-domain}. The shaded areas are the regions within three standard deviations of the sample mean, using the same $100$ draws as in Figure~\ref{fig:combustion-modes}. Training data is observed for $t \in [0.015, 0.017)$.}
    \label{fig:combustion-xvelocity}
\end{figure}

\begin{figure}[htbp]
    \centering
    \includegraphics[width=\textwidth]{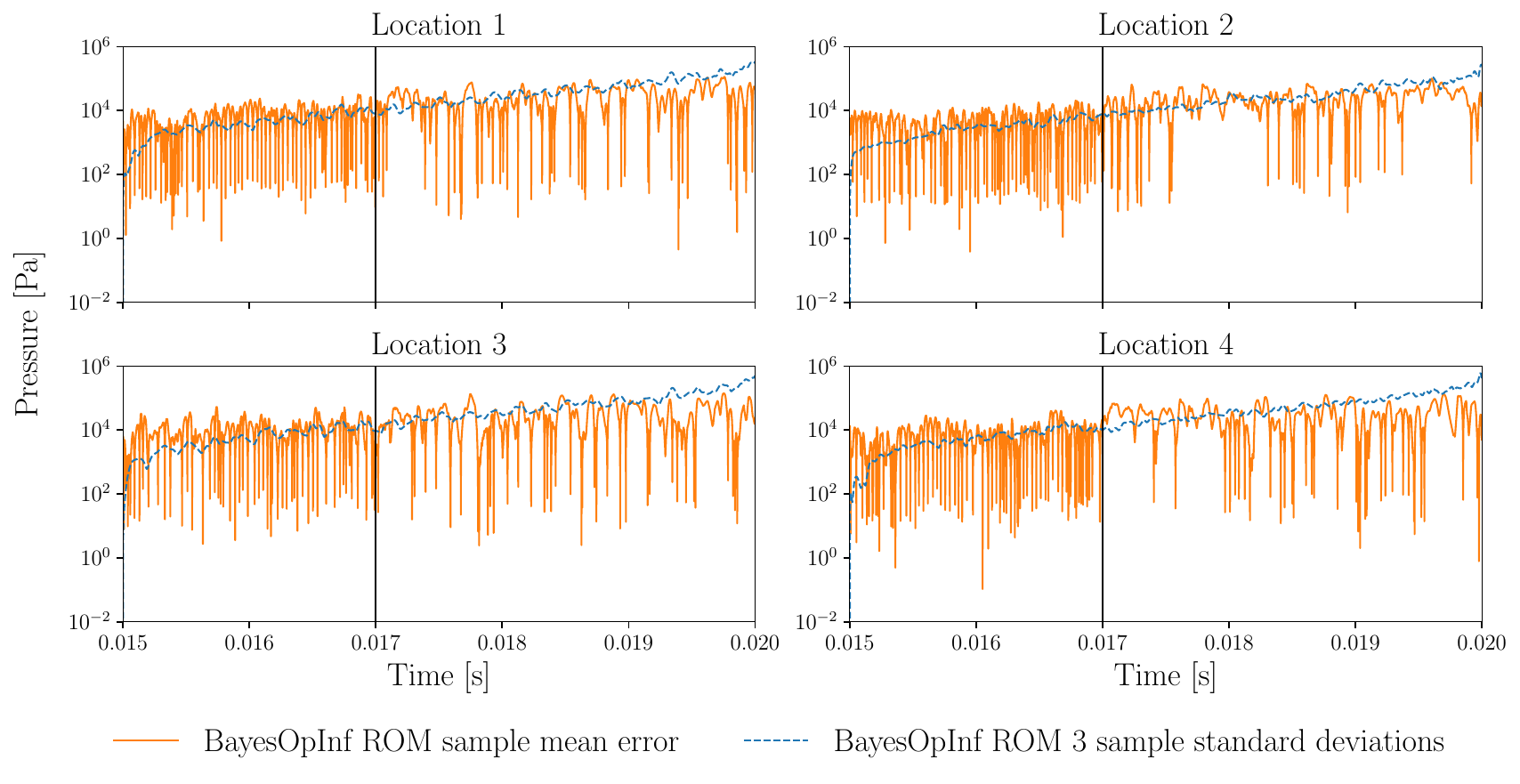}
    \vspace{-.75cm}
    \caption{Absolute error of the sample mean and three sample standard deviations for the pressure at the locations in Figure~\ref{fig:combustion-domain}, computed from $100$ draws from the posterior distribution of the Bayesian OpInf ROM. Note the logarithmic scale of the $y$-axis.}
    \label{fig:combustion-timeVspread-pressure}
    \bigskip\bigskip
    \includegraphics[width=\textwidth]{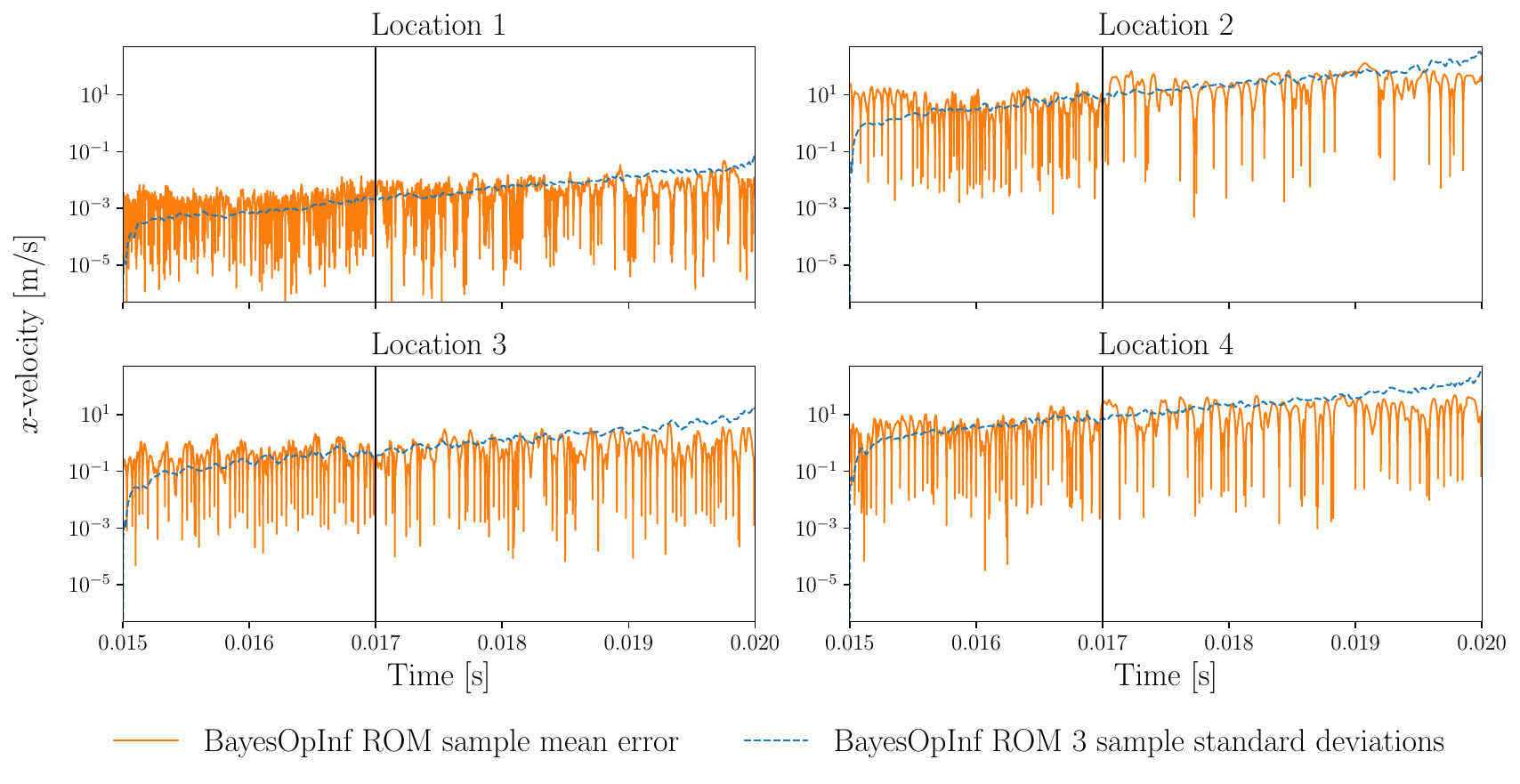}
    \vspace{-.75cm}
    \caption{Absolute error of the sample mean and three sample standard deviations for the horizontal velocity at the locations in Figure~\ref{fig:combustion-domain}, computed from $100$ draws from the posterior distribution of the Bayesian OpInf ROM. Note the logarithmic scale of the $y$-axis.}
    \label{fig:combustion-timeVspread-xvelocity}
\end{figure}

Figures \ref{fig:combustion-pressure} and \ref{fig:combustion-xvelocity} show the pressure and $x$-velocity, respectively, of the reconstructed ROM trajectory at the sensor locations, together with the data from the high-dimensional solver. We examine these variables in particular because they exhibit low projection error in the prediction regime with respect to the underlying POD basis $\mathbf{V}$. The full-order trajectory computed by GEMS is also shown for comparison. As in Figure~\ref{fig:combustion-modes}, the uncertainty bands remain quite tight in the training regime but start to noticeably increase in width shortly after the transition to the prediction regime. The bands are widest at the peaks and troughs of the pressure signal, reflecting that the model more accurately captures the frequency than the amplitude. While the sample mean is sometimes point-wise inaccurate, the true solution is generally well captured within three standard deviations of the sample mean.
Figures \ref{fig:combustion-timeVspread-pressure} and \ref{fig:combustion-timeVspread-xvelocity} show this more precisely for pressure and $x$-velocity, respectively, by plotting the absolute error
and three sample standard deviation (half the width of the uncertainty bands in Figures \ref{fig:combustion-pressure} and \ref{fig:combustion-xvelocity}) as a function of time at the sensor locations. Sharp drops in the sample mean error indicate that the sample mean crosses the true values. The width of the uncertainty band generally exceeds the absolute sample mean error in the prediction regime, especially for $t > 0.018$~s.
The sample mean error and the standard deviation are positively correlated in time, i.e., the uncertainty in the prediction generally increases with the error in the model.

\section{Concluding remarks}
\label{sec:conclusion}

The proposed BayesOpInf method non-intrusively propagates uncertainties from training data noise and model form error to learned low-dimensional operators of a ROM. Specifically, a linear Bayesian inference driven by projected state data and Gaussian priors leads to Gaussian posterior distributions for the reduced operators. Because constructing the posterior distributions for the reduced operators and integrating particular realizations of the ROM are relatively inexpensive operations, we can efficiently equip the resulting ROM solutions with uncertainty bands through Monte Carlo sampling.

BayesOpInf is a natural generalization of the deterministic OpInf framework. Since the inference problem is generic and carried out only using the solution snapshots and their time-derivative data, BayesOpInf is a promising tool for more specialized applications of OpInf, including localized methods \cite{rudy2022localopinf}, problems with detailed structure \cite{benner2022opInfIncompressible}, or parametric settings \cite{MKW2021popinf}. Furthermore, because the posterior mean of BayesOpInf is equivalent to Tikhonov-regularized OpInf, the approach is applicable for any Tikhonov-type regularization strategy. Hierarchical Bayes models may be a good candidate for formulating an implicit, non-Tikhonov regularization in the future. For now, we observe that the choice of the Tikhonov regularization remains an important issue for the stability of ROMs, but coupling the state-of-the-art error-based regularization selection of \cite{MHW2021regOpInfCombustion} with BayesOpInf leads to strong performance. Even with carefully chosen regularization, reduced models drawn from the posterior distribution may be unstable. Ensuring ROM stability---and, more generally, embedding desired model structure such as symmetry or spectral properties---through the posterior distribution remains an important topic for future work.

Our Bayesian approach has an inherent connection to Gaussian processes (GPs). A further discussion on this fact is provided in Appendix 1, wherein a GP surrogate model is used for the likelihood definition in the Bayesian inference. One open question is whether additional known information about the problem structure can be used to more strongly inform the Bayesian inference, for example by imposing a non-Gaussian distribution to model the reduced-order residual. Another potential extension of this work is utilizing the uncertainty indicators in the ROM solutions to define acquisition functions for active learning, for example adaptively selecting new training snapshots over time and parameter domains.

In summary, BayesOpInf quantifies the modeling uncertainties stemming from the data noise, the discretization errors of time derivatives, and the model misspecification error when the full-order model generating the snapshots does not match the assumed form of the ROM, as well as the uncertainties introduced by the use of Tikhonov regularization. As demonstrated in the numerical examples, BayesOpInf is effective in constructing accurate ROMs from noisy snapshot data and is capable of capturing complex physics while simultaneously providing meaningful credible intervals.

\section*{Acknowledgments}

This work was supported in part by the ARPA-E DIFFERENTIATE Program under award DE-AR0001208; the U.S. Department of Energy, National Nuclear Security Administration under award DE-NA0003969; the U.S. Department of Energy AEOLUS MMICC center under award DE-SC0019303; and the Air Force Center of Excellence on Multi-Fidelity Modeling of Rocket Combustor Dynamics (Air Force Office of Scientific Research), award FA9550-17-1-0195. The first author also acknowledges the financial support from Sectorplan B\`{e}ta (the Netherlands) under the focus area \emph{Mathematics of Computational Science}. The authors express their appreciation to Dr.~Nicole~Aretz for her valuable feedback.

\renewcommand{\theequation}{A.\arabic{equation}}
\setcounter{equation}{0}
\section*{Appendix 1: Alternative definition of the likelihood}

In \eqref{eq:leastsqr-bayes}, an independent Gaussian noise was considered for the modeling of residual errors. If one would like to take into account the correlation over time and/or construct a surrogate model of the ROM closure, a Gaussian process can be used to replace the noise term, i.e.,
\begin{equation}\label{eq:klinreg}
\frac{\dd}{\dd t} \hat{q}_i(t) = \vb{d}(\hat{\vb{q}}(t), \vb{u}(t))\trp\hat{\vb{o}}_i +\epsilon_i(t), \quad \epsilon_i(t)\sim \mathcal{GP}(0,\sigma_i^2\kappa(t,t'|\vb*{\ell}_i))\,,
\end{equation}
in which $\kappa(t,t'|\vb*{\ell}_i)$ is a given kernel function with hyperparameters $\vb*{\ell}_i$. In this case the likelihood is alternatively defined as
\begin{equation}
p(\vb{r}_i|~\vb{D},\hat{\vb{o}}_i,\sigma_i^2,\vb*{\ell})=\mathcal{N}(\vb{r}_i |~\vb{D}\hat{\vb{o}}_i,\sigma_i^2\vb{K}_i(\vb*{\ell}_i))\,,
\end{equation}
where the covariance matrix $\vb{K}_i(\vb*{\ell}_i)\in \mathbb{R}^{k\times k}$ is defined to be $\vb{K}_i(\vb*{\ell}_i) = \kappa(\mathcal{T},\mathcal{T}|\vb*{\ell}_i)$ and assumed to be positive definite. Here $\mathcal{T}$ denotes the set of $k$ time instances $\{t_0,t_1,\cdots,t_{k-1}\}$ corresponding to the training data. Conditioning on the training data $(\vb{D}, \vb{r}_i)$, the posterior distribution $p(\hat{\vb{o}}_i | ~\vb{D},\vb{r}_i,\sigma_i^2, \vb*{\lambda}_i,\vb*{\ell}_i)=\mathcal{N}(\hat{\vb{o}}_i|\vb*{\mu}_i,\vb*{\Sigma}_i)$ thus has new formulations of mean vector and covariance matrix as follows:
\begin{equation}
\begin{split}
\vb*{\Sigma}_i & = \sigma_i^2\left[ \mathrm{diag}(\vb*{\lambda}_i) +\vb{D}\trp\vb{K}^{-1}_i\vb{D} \right]^{-1} \\
\vb*{\mu}_i & =\left[ \mathrm{diag}(\vb*{\lambda}_i) +\vb{D}\trp\vb{K}^{-1}_i\vb{D} \right]^{-1} \left( \mathrm{diag}(\vb*{\lambda}_i)\vb*{\beta}_i+\vb{D}\trp\vb{K}^{-1}_i\vb{r}_i  \right)\\
& = \vb*{\beta}_i + \underbrace{\left[ \mathrm{diag}(\vb*{\lambda}_i) +\vb{D}\trp\vb{K}^{-1}_i\vb{D} \right]^{-1}\vb{D}\trp\vb{K}^{-1}_i(\vb{r}_i -\vb{D}\vb*{\beta}_i)}_{:=\delta\vb*{\mu}_i}\,.
\end{split}
\end{equation}
Note that the correction term $\delta \vb*{\mu}_i= \vb*{\mu}_i - \vb*{\beta}_i$ is then the solution to the following generalized least squares problem with a Tikhonov regularization
\begin{equation}
\delta \vb*{\mu}_i = \underset{\vb*{\eta}\in \mathbb{R}^{d(r,m)}}{\textrm{arg\,min}}\left\{ \|\vb{r}_i -\vb{D}(\vb*{\beta}_i +\vb*{\eta})\|_{\vb{K}^{-1}_i}^2  + \|\mathrm{diag}(\vb*{\lambda}_i)^{\frac{1}{2}}\vb*{\eta}\|_2^2   \right\}\,.
\end{equation}
Here the norm $\|\cdot\|_{\vb{K}^{-1}_i}$ is defined as $\|\vb{r}\|_{\vb{K}^{-1}_i}=\sqrt{\vb{r}\trp\vb{K}^{-1}_i\vb{r}}$, $\vb{r}\in\mathbb{R}^k$. In the meantime, the estimate of hyperparameters $(\sigma_i^2,\vb*{\lambda}_i,\vb*{\ell}_i)$ by maximizing marginal likelihood is reformulated as
\begin{equation}
\begin{split}
(\sigma^{*2}_i,\vb*{\lambda}_i^*,\vb*{\ell}^*_i)  =  \underset{\sigma_i^2,\vb*{\lambda}_i,\vb*{\ell}_i}{\textrm{arg\,max}}&~\log~p(\vb{r}_i|~\vb{D},\sigma_i^2,\vb*{\lambda}_i,\vb*{\ell}_i) \\
= \underset{\sigma_i^2,\vb*{\lambda}_i,\vb*{\ell}_i}{\textrm{arg\,max}}&~\log\int p(\vb{r}_i|~\vb{D},\hat{\vb{o}}_i,\sigma_i^2,\vb*{\ell}_i)~ p(\hat{\vb{o}}_i|~\sigma_i^2,\vb*{\lambda}_i)~\dd\hat{\vb{o}}_i \\
 = \underset{\sigma_i^2,\vb*{\lambda}_i,\vb*{\ell}_i}{\textrm{arg\,max}}&~\log~\mathcal{N}\left(\vb{r}_i | ~\vb{D}\vb*{\beta}_i, \sigma_i^2\left(\vb{D}\mathrm{diag}(\vb*{\lambda}_i)^{-1}\vb{D}\trp+\vb{K}_i(\vb*{\ell}_i)\right)\right)\,.
\end{split}
\end{equation}

The model~\eqref{eq:klinreg} is equivalent to a universal kriging with the following  Gaussian process prior:
\begin{equation}
\frac{\dd}{\dd t}  \hat{q}_i(t) \sim \mathcal{GP}\bigg(\vb{d}(\hat{\vb{q}}(t), \vb{u}(t))\trp\vb*{\beta}_i, \sigma_i^2\vb{d}(\hat{\vb{q}}(t), \vb{u}(t))\trp\mathrm{diag}(\vb*{\lambda}_i)^{-1}\vb{d}(\hat{\vb{q}}(t'), \vb{u}(t'))+\sigma_i^2\kappa(t,t'|\vb*{\ell}_i)\bigg)\,.
\end{equation}
Accordingly, the posterior prediction of $\frac{\dd}{\dd t}  \hat{q}_i(t)$ is again a Gaussian process whose mean and covariance functions take the form of
\begin{equation}
\begin{split}
\mathbb{E}\left[\frac{\dd}{\dd t}  \hat{q}_i(t) \middle\vert~\vb{D},\vb{r}_i\right] = &~ \vb{d}(\hat{\vb{q}}(t), \vb{u}(t))\trp\vb*{\mu}_i+\kappa(t,\mathcal{T}|\vb*{\ell}_i)\vb{K}^{-1}_i(\vb{r}_i-\vb{D}\vb*{\mu}_i)\,,\quad \text{and}\\
\mathbb{C}ov\left[\frac{\dd}{\dd t}  \hat{q}_i(t),\frac{\dd}{\dd t}  \hat{q}_i(t') \middle\vert~\vb{D},\vb{r}_i\right]  = &~ \sigma_i^2\kappa(t,t'|\vb*{\ell}_i)-\sigma_i^2\kappa(t,\mathcal{T}|\vb*{\ell}_i)\vb{K}^{-1}_i\kappa(\mathcal{T},t'|\vb*{\ell}_i)\\
+\big[\vb{d}(\hat{\vb{q}}(t), \vb{u}(t))\trp&-\kappa(t,\mathcal{T}|\vb*{\ell}_i)\vb{K}_i^{-1}\vb{D}\big]\vb*{\Sigma}_i\big[\vb{d}(\hat{\vb{q}}(t'), \vb{u}(t'))-\vb{D}\trp\vb{K}_i^{-1}\kappa(\mathcal{T},t'|\vb*{\ell}_i)\big]\,,
\end{split}
\end{equation}
respectively. Note that the second term in the mean function, i.e., $\kappa(t,\mathcal{T}|\vb*{\ell}_i)\vb{K}^{-1}_i(\vb{r}_i-\vb{D}\vb*{\mu}_i)$, formulates a surrogate of the closure in the linear model~\eqref{eq:reducedorder}.

\renewcommand{\theequation}{B.\arabic{equation}}
\setcounter{equation}{0}
\section*{Appendix 2: Error-based regularization selection}

Algorithm~\ref{alg:OpInf-Reg} details the OpInf procedure where the regularization parameters are selected with an optimization routine \cite{MHW2021regOpInfCombustion}.
Regularization parameters candidates are assessed by solving the corresponding deterministic OpInf problem~\eqref{eq:leastsqr}, then using the resultant ROM to reconstruct the training data for $L \ge 1$ initial conditions and make predictions in time beyond the training regime.
The projected training snapshots are organized by initial condition, i.e., $\hat{\mathbf{Q}}^{(\ell)}$ contains all projected snapshots for the $\ell$th initial condition.

\begin{algorithm}[t]
\begin{algorithmic}[1]
\Procedure{OpInf}{
    Projected training snapshots $\hat{\mathbf{Q}} = [~\hat{\mathbf{Q}}^{(1)}~\cdots~\hat{\mathbf{Q}}^{(L)}~]\in\mathbb{R}^{r\times k}$,
    training inputs $\mathbf{U}\in\mathbb{R}^{m\times k}$,
    projected snapshot time derivatives $\mathbf{R} = [~\mathbf{r}_{1}~\cdots~\mathbf{r}_{r}~]\trp\in\mathbb{R}^{r\times k}$,
    final time $t_{f} > t_{0}$,
    bound margin $\tau \ge 1$
    }
    \vspace{.0625in}
    \State $\mathbf{D} \gets \texttt{data\_matrix}(\hat{\vb{Q}}, \mathbf{U})$
        \Comment{Construct data matrix as in \eqref{eq:matrix-definitions}.}
        \label{step:data-matrix}
    \State $B \gets \tau \max_{i,j,\ell}|\hat{\mathbf{Q}}_{i,j}^{(\ell)}|$
        \Comment{Set a bound for the ROM stability criterion.}

    \vspace{.0625in}
    \LineComment{Subprocedure for assessing regularization parameters.}
    \Procedure{OpInfError}{
    regularizer $\boldsymbol{\lambda}\in\mathbb{R}^{d(r,m)}$
    }
    \For{$i= 1, \ldots, r$}
        \Comment{OpInf with current regularization.}
        \State $\hat{\mathbf{o}}_{i} \gets \arg\min_{\boldsymbol{\mu}_{i}}\left\|\mathbf{D}\boldsymbol{\mu}_{i} - \mathbf{r}_{i}\right\|_{2}^{2} + \left\|\textrm{diag}(\boldsymbol{\lambda})^{1/2}\boldsymbol{\mu}_{i}\right\|_{2}^{2}$
    \EndFor
    \For{$\ell = 1, \ldots, L$}
        \State $\tilde{\mathbf{Q}} \gets~\textrm{integrate \eqref{eq:reducedorder} with}~\hat{\mathbf{O}} = [~\hat{\mathbf{o}}_{1}~\cdots~\hat{\mathbf{o}}_{r}~]\trp$ from $\hat{\mathbf{q}}_{0} = \hat{\mathbf{Q}}_{:,0}$ over $[t_0,t_f]$
        \label{step:integrate-rom}
    \EndFor
    \If{$\max_{i,j,\ell}|\tilde{\mathbf{Q}}^{(\ell)}_{i,j}| > B$}
        \Comment{Ensure the ROM is stable on $[t_0,t_f]$ for all ICs.}
        \State \textbf{return}~$\infty$
    \EndIf
    \State \textbf{return}~$\displaystyle\frac{1}{L}\sum_{\ell=1}^{L}\left\|\hat{\mathbf{Q}}^{(\ell)} - \tilde{\mathbf{Q}}_{:,:k}^{(\ell)}\right\|$
        \Comment{Compare ROM outputs to training data.}
    \EndProcedure

    \LineComment{Minimize the subprocedure to choose the regularization parameters.}
    \State \textbf{return}~$\arg\min_{\boldsymbol{\lambda}}$~\textproc{OpInfError}$(\boldsymbol{\lambda})$
    \label{step:minimize-error}
\EndProcedure
\end{algorithmic}
\caption{Operator Inference with error-based regularization selection (Algorithm~1 of~\cite{MHW2021regOpInfCombustion}).}
\label{alg:OpInf-Reg}
\end{algorithm}

In the results for Section~\ref{sec:euler} (Figure~\ref{fig:euler-basisVerror}), Algorithm~\ref{alg:OpInf-Reg} is tailored to the purely quadratic ROM~\eqref{eq:rom-quadratic}. In this setting, there are no external inputs $\mathbf{U}$, the data matrix of step~\ref{step:data-matrix} is given by $\mathbf{D} = (\hat{\mathbf{Q}}\odot\hat{\mathbf{Q}})\trp$, and the operator matrix of step~\ref{step:integrate-rom} is $\hat{\mathbf{O}} = \hat{\mathbf{H}}$.
The regularization parameters $\boldsymbol{\lambda}$ are globally parameterized by a single scalar hyperparameter $\lambda \ge 0$, i.e., $\boldsymbol{\lambda} = [~\lambda~\cdots~\lambda~]\trp$. Hence, the minimization in step~\ref{step:minimize-error} can be carried out with a one-dimensional optimization on $\lambda$ via, e.g., Brent's method~\cite{Brent2002minimization}.

For the combustion problem of Section~\ref{sec:combustion}, we have a single initial condition ($L = 1$) and the ROM is \eqref{eq:reducedorder} with data matrix $\mathbf{D} = [~\hat{\vb{Q}}\trp~~(\hat{\vb{Q}}\odot \hat{\vb{Q}})\trp~~\vb{U}\trp~~\vb{1}_k~]$ and reduced-order operators $\hat{\mathbf{O}} = [~\hat{\mathbf{A}}~\hat{\mathbf{H}}~\hat{\mathbf{B}}~\hat{\mathbf{c}}~]$.
The regularization parameters are parameterized by two scalars $\lambda_{1},\lambda_{2}\ge 0$:
\begin{equation}
    \boldsymbol{\lambda}
    = [~\underbrace{\lambda_{1},\ldots,\lambda_{1}}_{r~\textrm{times}},\underbrace{\lambda_{2},\ldots,\lambda_{2}}_{\frac{r(r+1)}{2}~\textrm{times}},\lambda_{1},\lambda_{1}~]\trp.
\end{equation}
This choice of $\boldsymbol{\lambda}$ penalizes the entries of the quadratic operator $\hat{\mathbf{H}}$ by $\lambda_{2}$ and the entries of the remaining operators $\hat{\mathbf{A}}$, $\hat{\mathbf{B}}$, and $\hat{\mathbf{c}}$ by $\lambda_{1}$. More precisely, the regularization term as written in \eqref{eq:opinf} is given by
\begin{equation}
    \mathcal{P}\left([~\hat{\vb{A}}~\hat{\vb{H}}~\hat{\vb{B}}~\hat{\vb{c}}~]\right)
    = \lambda_{1}\left(
        \|\hat{\mathbf{A}}\|_{F}^{2}
        + \|\hat{\mathbf{B}}\|_{F}^{2}
        + \|\hat{\mathbf{c}}\|_{2}^{2}
    \right) + \lambda_{2}\|\hat{\mathbf{H}}\|_{F}^{2}.
\end{equation}
The optimization in step~\ref{step:minimize-error} of Algorithm~\ref{alg:OpInf-Reg} is carried out with a coarse grid search in the principal quadrant of $\mathbb{R}^{2}$, followed by the Nelder-Mead method~\cite{nelder1965simplex}.

\bibliographystyle{abbrv}
\bibliography{refs.bib}

\begin{thebibliography}{10}

\bibitem{oden2011fembook}
M.~Ainsworth and J.~T. Oden.
\newblock {\em A Posteriori Error Estimation in Finite Element Analysis},
  volume~37.
\newblock John Wiley \& Sons, 2011.

\bibitem{oden2004verification}
I.~Babuska and J.~T. Oden.
\newblock Verification and validation in computational engineering and science:
  basic concepts.
\newblock {\em Computer Methods in Applied Mechanics and Engineering},
  36(193):4057--4066, 2004.

\bibitem{benner2022opInfIncompressible}
P.~Benner, P.~Goyal, J.~Heiland, and I.~P. Duff.
\newblock Operator inference and physics-informed learning of low-dimensional
  models for incompressible flows.
\newblock {\em Electronic Transactions on Numerical Analysis}, 56:28--51, 2022.

\bibitem{BGW2015pmorSurvery}
P.~Benner, S.~Gugercin, and K.~Willcox.
\newblock A survey of projection-based model reduction methods for parametric
  dynamical systems.
\newblock {\em SIAM Review}, 57(4):483--531, 2015.

\bibitem{Berkooz1993}
G.~Berkooz, P.~Holmes, and J.~Lumley.
\newblock The proper orthogonal decomposition in the analysis of turbulent
  flows.
\newblock {\em Annual Review of Fluid Mechanics}, 25:539--575, 1993.

\bibitem{box2011bayesian}
G.~E.~P. Box and G.~C. Tiao.
\newblock {\em Bayesian Inference in Statistical Analysis}.
\newblock John Wiley \& Sons, 2011.

\bibitem{Brent2002minimization}
R.~P. Brent.
\newblock {\em Algorithms for Minimization Without Derivatives}.
\newblock Dover Publications, New York, NY, 2002.

\bibitem{champion2019sindy}
K.~Champion, B.~Lusch, J.~N. Kutz, and S.~L. Brunton.
\newblock Data-driven discovery of coordinates and governing equations.
\newblock {\em Proceedings of the National Academy of Sciences},
  116(45):22445--22451, 2019.

\bibitem{chkrebtii2016bayesian}
O.~A. Chkrebtii, D.~A. Campbell, B.~Calderhead, and M.~A. Girolami.
\newblock Bayesian solution uncertainty quantification for differential
  equations.
\newblock {\em Bayesian Analysis}, 11(4):1239--1267, 2016.

\bibitem{DeBrabanter2013localddt}
K.~De~Brabanter, J.~De~Brabanter, B.~De~Moor, and I.~Gijbels.
\newblock Derivative estimation with local polynomial fitting.
\newblock {\em Journal of Machine Learning Research}, 14(1), 2013.

\bibitem{del2008proper}
D.~del Castillo-Negrete, D.~A. Spong, and S.~P. Hirshman.
\newblock Proper orthogonal decomposition methods for noise reduction in
  particle-based transport calculations.
\newblock {\em Physics of Plasmas}, 15(9):092308, 2008.

\bibitem{epps2010PODthresh}
B.~P. Epps and A.~H. Techet.
\newblock An error threshold criterion for singular value decomposition modes
  extracted from {PIV} data.
\newblock {\em Experiments in fluids}, 48(2):355--367, 2010.

\bibitem{oden2015bayesatomistic}
K.~Farrell, J.~T. Oden, and D.~Faghihi.
\newblock A {B}ayesian framework for adaptive selection, calibration, and
  validation of coarse-grained models of atomistic systems.
\newblock {\em Journal of Computational Physics}, 295:189--208, 2015.

\bibitem{rudy2022localopinf}
R.~Geelen and K.~Willcox.
\newblock Localized non-intrusive reduced-order modelling in the operator
  inference framework.
\newblock {\em Philosophical Transactions of the Royal Society A: Mathematical,
  Physical and Engineering Sciences}, 380(2229):20210206, 2022.

\bibitem{GW2021learning}
O.~Ghattas and K.~Willcox.
\newblock Learning physics-based models from data: perspectives from inverse
  problems and model reduction.
\newblock {\em Acta Numerica}, 30:445--554, 2021.

\bibitem{GPT1999vortexsheddingPOD}
W.~R. Graham, J.~Peraire, and K.~Y. Tang.
\newblock Optimal control of vortex shedding using low-order models. {P}art
  {I}---{O}pen-loop model development.
\newblock {\em International Journal for Numerical Methods in Engineering},
  44(7):945--972, 1999.

\bibitem{guo2019data}
M.~Guo and J.~S. Hesthaven.
\newblock Data-driven reduced order modeling for time-dependent problems.
\newblock {\em Computer Methods in Applied Mechanics and Engineering},
  345:75--99, 2019.

\bibitem{HMPT2011rNLA}
N.~Halko, P.-G. Martinsson, and J.~A. Tropp.
\newblock Finding structure with randomness: {P}robabilistic algorithms for
  constructing approximate matrix decompositions.
\newblock {\em SIAM Review}, 53(2):217--288, 2011.

\bibitem{HHSFAMT2015gems}
M.~E. Harvazinski, C.~Huang, V.~Sankaran, T.~W. Feldman, W.~E. Anderson, C.~L.
  Merkle, and D.~G. Talley.
\newblock Coupling between hydrodynamics, acoustics, and heat release in a
  self-excited unstable combustor.
\newblock {\em Physics of Fluids}, 27(4):045102, 2015.

\bibitem{oden2013bayestumor}
A.~Hawkins-Daarud, S.~Prudhomme, K.~G. van~der Zee, and J.~T. Oden.
\newblock Bayesian calibration, validation, and uncertainty quantification of
  diffuse interface models of tumor growth.
\newblock {\em Journal of Mathematical Biology}, 67(6):1457--1485, 2013.

\bibitem{hesthaven2016certified}
J.~S. Hesthaven, G.~Rozza, and B.~Stamm.
\newblock {\em Certified reduced basis methods for parametrized partial
  differential equations}, volume 590.
\newblock Springer, 2016.

\bibitem{HBDK2021uqsindy}
S.~M. Hirsh, D.~A. Barajas-Solano, and J.~N. Kutz.
\newblock Sparsifying priors for {B}ayesian uncertainty quantification in model
  discovery.
\newblock {\em arXiv preprint arXiv:2107.02107}, 2021.

\bibitem{HDKM2018combustion-is-hard}
C.~Huang, K.~Duraisamy, and C.~Merkle.
\newblock Challenges in reduced order modeling of reacting flows.
\newblock In {\em 2018 Joint Propulsion Conference}, Cincinnati, OH, 2018.
\newblock Paper AIAA-2018-4675.

\bibitem{HDM2019poddeim-robustness}
C.~Huang, K.~Duraisamy, and C.~L. Merkle.
\newblock Investigations and improvement of robustness of reduced-order models
  of reacting flow.
\newblock {\em AIAA Journal}, 57(12):5377--5389, 2019.

\bibitem{HXDM2018rocketrom-poddeim}
C.~Huang, J.~Xu, K.~Duraisamy, and C.~Merkle.
\newblock Exploration of reduced-order models for rocket combustion
  applications.
\newblock In {\em 2018 AIAA Aerospace Sciences Meeting}, Orlando, FL, 2018.
\newblock Paper AIAA-2018-1183.

\bibitem{jain2021performance}
P.~Jain, S.~McQuarrie, and B.~Kramer.
\newblock Performance comparison of data-driven reduced models for a
  single-injector combustion process.
\newblock In {\em AIAA Propulsion and Energy 2021 Forum}, Virtual Event, 2021.
\newblock Paper AIAA-2021-3633.

\bibitem{oden2020bayescovid}
P.~K. Jha, L.~Cao, and J.~T. Oden.
\newblock Bayesian-based predictions of {COVID}-19 evolution in {T}exas using
  multispecies mixture-theoretic continuum models.
\newblock {\em Computational Mechanics}, 66(5):1055--1068, 2020.

\bibitem{KhatriRoa1968product}
C.~Khatri and C.~R. Rao.
\newblock Solutions to some functional equations and their applications to
  characterization of probability distributions.
\newblock {\em Sankhy{\=a}: The Indian Journal of Statistics, Series A}, pages
  167--180, 1968.

\bibitem{kolda2009tensor}
T.~G. Kolda and B.~W. Bader.
\newblock Tensor decompositions and applications.
\newblock {\em SIAM review}, 51(3):455--500, 2009.

\bibitem{Lee2020autoencoder}
K.~Lee and K.~T. Carlberg.
\newblock Model reduction of dynamical systems on nonlinear manifolds using
  deep convolutional autoencoders.
\newblock {\em Journal of Computational Physics}, 404:108973, 2020.

\bibitem{MHW2021regOpInfCombustion}
S.~A. McQuarrie, C.~Huang, and K.~E. Willcox.
\newblock Data-driven reduced-order models via regularised operator inference
  for a single-injector combustion process.
\newblock {\em Journal of the Royal Society of New Zealand}, 51(2):194--211,
  2021.

\bibitem{MKW2021popinf}
S.~A. McQuarrie, P.~Khodabakhshi, and K.~E. Willcox.
\newblock Non-intrusive reduced-order models for parametric partial
  differential equations via data-driven operator inference.
\newblock {O}den {I}nstitute {R}eport 21-17, University of Texas at Austin,
  2021.

\bibitem{nelder1965simplex}
J.~A. Nelder and R.~Mead.
\newblock A simplex method for function minimization.
\newblock {\em The Computer Journal}, 7(4):308--313, 1965.

\bibitem{oden2018predictivemodelling}
J.~T. Oden.
\newblock Adaptive multiscale predictive modelling.
\newblock {\em Acta Numerica}, 27:353--450, 2018.

\bibitem{oden2017predictive}
J.~T. Oden, I.~Babu{\v{s}}ka, and D.~Faghihi.
\newblock Predictive computational science: Computer predictions in the
  presence of uncertainty.
\newblock {\em Encyclopedia of Computational Mechanics Second Edition}, pages
  1--26, 2017.

\bibitem{oden2005uncertainty}
J.~T. Oden, I.~Babu{\v{s}}ka, F.~Nobile, Y.~Feng, and R.~Tempone.
\newblock Theory and methodology for estimation and control of errors due to
  modeling, approximation, and uncertainty.
\newblock {\em Computer Methods in Applied Mechanics and Engineering},
  194(2-5):195--204, 2005.

\bibitem{oden2002estimation}
J.~T. Oden and S.~Prudhomme.
\newblock Estimation of modeling error in computational mechanics.
\newblock {\em Journal of Computational Physics}, 182(2):496--515, 2002.

\bibitem{pan2015sparse}
W.~Pan, Y.~Yuan, J.~Gon{\c{c}}alves, and G.-B. Stan.
\newblock A sparse {B}ayesian approach to the identification of nonlinear
  state-space systems.
\newblock {\em IEEE Transactions on Automatic Control}, 61(1):182--187, 2015.

\bibitem{Peherstorfer2020reprojection}
B.~Peherstorfer.
\newblock Sampling low-dimensional {M}arkovian dynamics for pre-asymptotically
  recovering reduced models from data with operator inference.
\newblock {\em SIAM Journal on Scientific Computing}, 42(5):A3489--A3515, 2020.

\bibitem{PW2016operatorInference}
B.~Peherstorfer and K.~Willcox.
\newblock Data-driven operator inference for nonintrusive projection-based
  model reduction.
\newblock {\em Computer Methods in Applied Mechanics and Engineering},
  306:196--215, 2016.

\bibitem{oden2015bayesdamage}
E.~Prudencio, P.~Bauman, D.~Faghihi, K.~Ravi-Chandar, and J.~T. Oden.
\newblock A computational framework for dynamic data-driven material damage
  control, based on {B}ayesian inference and model selection.
\newblock {\em International Journal for Numerical Methods in Engineering},
  102(3-4):379--403, 2015.

\bibitem{QKPW2020liftAndLearn}
E.~Qian, B.~Kramer, B.~Peherstorfer, and K.~Willcox.
\newblock Lift \& {L}earn: {P}hysics-informed machine learning for large-scale
  nonlinear dynamical systems.
\newblock {\em Physica {D}: {N}onlinear {P}henomena}, 406:132401, 2020.

\bibitem{quarteroni2014reduced}
A.~Quarteroni, G.~Rozza, et~al.
\newblock {\em Reduced order methods for modeling and computational reduction},
  volume~9.
\newblock Springer, 2014.

\bibitem{schaeffer2017learning}
H.~Schaeffer.
\newblock Learning partial differential equations via data discovery and sparse
  optimization.
\newblock {\em Proceedings of the Royal Society A: Mathematical, Physical and
  Engineering Sciences}, 473(2197):20160446, 2017.

\bibitem{schmid2010dynamic}
P.~J. Schmid.
\newblock Dynamic mode decomposition of numerical and experimental data.
\newblock {\em Journal of Fluid Mechanics}, 656:5--28, 2010.

\bibitem{sirovich1987turbulence}
L.~Sirovich.
\newblock Turbulence and the dynamics of coherent structures. {I}. {C}oherent
  structures.
\newblock {\em Quarterly of Applied Mathematics}, 45(3):561--571, 1987.

\bibitem{smith2013uncertainty}
R.~C. Smith.
\newblock {\em Uncertainty Quantification: Theory, Implementation, and
  Applications}.
\newblock SIAM, 2013.

\bibitem{SKHW2020romCombustion}
R.~Swischuk, B.~Kramer, C.~Huang, and K.~Willcox.
\newblock Learning physics-based reduced-order models for a single-injector
  combustion process.
\newblock {\em AIAA Journal}, 58(6):2658--2672, 2020.

\bibitem{uy2021probabilistic}
W.~I.~T. Uy and B.~Peherstorfer.
\newblock Probabilistic error estimation for non-intrusive reduced models
  learned from data of systems governed by linear parabolic partial
  differential equations.
\newblock {\em ESAIM: Mathematical Modelling and Numerical Analysis},
  55(3):735--761, 2021.

\bibitem{uy2021activeOpInf}
W.~I.~T. Uy, Y.~Wang, Y.~Wen, and B.~Peherstorfer.
\newblock Active operator inference for learning low-dimensional
  dynamical-system models from noisy data.
\newblock {\em arXiv preprint arXiv:2107.09256}, 2021.

\bibitem{vanLoan2000kronecker}
C.~F. Van~Loan.
\newblock The ubiquitous {K}ronecker product.
\newblock {\em Journal of Computational and Applied Mathematics},
  123(1-2):85--100, 2000.

\bibitem{venturi2006PODperturb}
D.~Venturi.
\newblock On proper orthogonal decomposition of randomly perturbed fields with
  applications to flow past a cylinder and natural convection over a horizontal
  plate.
\newblock {\em Journal of Fluid Mechanics}, 559:215--254, 2006.

\bibitem{WD1981oxidation}
C.~K. Westbrook and F.~L. Dryer.
\newblock Simplified reaction mechanisms for the oxidation of hydrocarbon fuels
  in flames.
\newblock {\em Combustion Science and Technology}, 27(1-2):31--43, 1981.

\bibitem{rasmussen2006gaussian}
C.~K. Williams and C.~E. Rasmussen.
\newblock {\em Gaussian Processes for Machine Learning}.
\newblock MIT Press, 2006.

\bibitem{zhang2018robust}
S.~Zhang and G.~Lin.
\newblock Robust data-driven discovery of governing physical laws with error
  bars.
\newblock {\em Proceedings of the Royal Society A: Mathematical, Physical and
  Engineering Sciences}, 474(2217):20180305, 2018.

\bibitem{zhuang2021model}
Q.~Zhuang, J.~M. Lorenzi, H.-J. Bungartz, and D.~Hartmann.
\newblock Model order reduction based on {R}unge-{K}utta neural networks.
\newblock {\em Data-Centric Engineering}, 2, 2021.

\end{thebibliography}

\end{document}